\documentclass[12pt]{amsart}
\usepackage[pdftex]{hyperref}
\hypersetup{colorlinks,citecolor=blue,linktocpage}
\usepackage{amsmath,amscd,amssymb,euscript}
\usepackage{mathtools}
\usepackage{latexsym}
\usepackage{graphicx}
\usepackage{tikz}
\usetikzlibrary{backgrounds}
\usetikzlibrary{decorations.fractals}
\usetikzlibrary{calc,intersections,through,backgrounds}
\usetikzlibrary{arrows.meta}
\tikzset{In/.tip = {Hooks[right]}}
\tikzset{Onto/.tip = {To[sep] To}}
\tikzset{Eq/.style = {-, double equal sign distance}}
\tikzset{Iso/.style = {edge node = {node[above] {$\sim$}}, inner sep = 0}}
\usetikzlibrary{positioning}
\usepackage{mathrsfs}
\usepackage{epstopdf}
\usepackage{rotating}
\usepackage{lscape}
\usepackage[all,cmtip]{xy}

\newcommand{\del}{\partial}
\newcommand{\bC}{{\mathbb C}}

\newcommand{\bH}{{\mathbb H}}

\newcommand{\bN}{{\mathbb N}}
\newcommand{\bO}{{\mathbb O}}
\newcommand{\bP}{{\mathbb P}}
\newcommand{\bQ}{{\mathbb Q}}
\newcommand{\bR}{{\mathbb R}}

\newcommand{\bZ}{{\mathbb Z}}

\newcommand{\cF}{{\mathcal F}}

\newcommand{\cM}{{\mathcal M}}

\newcommand{\cO}{{\mathcal O}}

\newcommand{\cX}{{\mathcal X}}

\newcommand{\Fgl}{\mathfrak{gl}}

\newcommand{\Fhol}{\mathfrak{hol}}
\newcommand{\FS}{\mathfrak{S}}
\newcommand{\oa}{{\overline{a}}}
\newcommand{\oA}{{\overline{A}}}
\newcommand{\oalpha}{{\overline{\alpha}}}

\newcommand{\ocJ}{{\overline{\mathcal J}}}

\newcommand{\oQ}{{\overline{Q}}}

\newcommand{\osigma}{{\overline{\sigma}}}

\newcommand{\tq}{{\widetilde{q}}}

\newcommand{\ra}{\rightarrow}
\newcommand{\lra}{\longrightarrow}

\newcommand{\surj}{\mathrel{\mathrlap{\rightarrow}\mkern1mu\rightarrow}}

\newcommand{\hint}[1]{\textcolor{gray!60!}{\textbf{ }}}

\def\cl{{\colon}}

\def\hC{{\widehat{C}}}

\newcommand{\Def}{\operatorname{Def}}
\newcommand{\Diff}{\operatorname{Diff}}
\newcommand{\Div}{\operatorname{Div}}
\newcommand{\End}{\operatorname{End}}

\newcommand{\GL}{\operatorname{GL}}
\newcommand{\Hol}{\operatorname{Hol}}
\newcommand{\Hom}{\operatorname{Hom}}

\newcommand{\Img}{\operatorname{Im}}
\newcommand{\Isom}{\operatorname{Isom}}
\newcommand{\Ker}{\operatorname{Ker}}

\newcommand{\Pic}{\operatorname{Pic}}

\newcommand{\rank}{\operatorname{rank}}
\newcommand{\Ric}{\operatorname{Ric}}
\newcommand{\slope}{\operatorname{slope}}

\newcommand{\Sym}{\operatorname{Sym}}
\newcommand{\td}{\operatorname{td}}
\newcommand{\Teich}{\operatorname{Teich}}
\newcommand{\trace}{\operatorname{tr}}
\newcommand{\vol}{\operatorname{vol}}

\usepackage[nameinlink,capitalize]{cleveref}
\usepackage{enumitem}

\newlist{anumerate}{enumerate}{1}
\setlist[anumerate]{label=(\alph*)}
\crefname{subsection}{Paragraph}{Paragraphs}
\Crefname{subsection}{Paragraph}{Paragraphs}
\crefname{equation}{}{}
\Crefname{equation}{}{}
\crefname{anumeratei}{}{}
\Crefname{anumeratei}{}{}

\newcommand\eq\coloneqq
\newcommand\into\hookrightarrow

\newcommand\iso\cong
\newcommand\minus\setminus
\newcommand\onto\twoheadrightarrow

\newcommand\pr{{\mathrm{pr}}}

\newcommand\tensor\otimes

\renewcommand\tilde\widetilde

\newcommand{\tR}{{\widetilde{R}}}

\DeclareMathOperator\ch{ch}

\renewcommand\emptyset\varnothing

\theoremstyle{definition}

\newtheorem{proposition}{Proposition}[section]
\newtheorem{lemma}[proposition]{Lemma}

\newtheorem{theorem}[proposition]{Theorem}
\newtheorem{theoremi}{Theorem}
\newtheorem{conjecture}[theoremi]{Conjecture}
\newtheorem{definition}[proposition]{Definition}
\newtheorem{defprop}[proposition]{Definition and Proposition}
\newtheorem{corollary}[proposition]{Corollary}
\newtheorem{question}[proposition]{Question}

\newtheorem{problem}[proposition]{Problem}
\newtheorem{further}[proposition]{Further}

\numberwithin{equation}{section}

\begin{document}

\title[Hyperk\"ahler manifolds]{Hyperk\"ahler manifolds}

\author{Elham Izadi}

\address{Department of Mathematics, University of California San Diego, 9500 Gilman Drive \# 0112, La Jolla, CA 92093-0112, USA}

\email{eizadi@ucsd.edu}

\author{\\Exercises by: Samir Canning}

\address{Department of Mathematics, University of California San Diego, 9500 Gilman Drive \# 0112, La Jolla, CA 92093-0112, USA}

\email{srcannin@ucsd.edu}

\author{Yajnaseni Dutta}

\address{Mathematik Zentrum\\Universit\"at Bonn,
Endenicher Allee 60, Germany.}

\email{ydutta@uni-bonn.de}

\author{David Stapleton}

\address{Department of Mathematics, University of Michigan, 530 Church St, Ann Arbor, MI 48109
}

\email{dajost@umich.edu}

\thanks{}


\begin{abstract}

We give an elementary introduction to hyperk\"ahler manifolds, survey some of their interesting properties and some open problems.

\end{abstract}

\maketitle

\tableofcontents

\section*{Introduction}

The cohomology of a compact K\"ahler manifold has remarkable properties, abstractified in the modern notion of a (polarized) Hodge structure. While the datum of a Hodge structure of weight $1$ is equivalent to the datum of a compact complex torus, this is no longer the case in higher weights. In weight $2$ there are remarkable examples of compact K\"ahler manifolds which are, mostly, determined by the polarized Hodge structure on their second cohomology. These are the hyperk\"ahler manifolds: higher dimensional analogues of K3 surfaces. In these lecture notes, we give an elementary introduction to hyperk\"ahler manifolds and survey some of their interesting properties.

We start by reviewing the notions of tensors, connections, the curvature tensor, Ricci curvature and some of their properties. We define parallel transport, holonomy and the Levi-Civita connection. We also describe the constraints posed by the holonomy on the curvature tensor. We define (locally) symmetric spaces and state the main structure theorem for them. We then state De Rham's decomposition theorem for simply connected complete Riemannian manifolds and Berger's classification of the holonomy groups of nonsymmetric, complete, connected, irreducible Riemannian manifolds. Berger's classification shows that hyperk\"ahler manifolds are the nonsymmetric complete connected irreducible Riemannian manifolds with holonomy group contained in $Sp(r)$: the group of automorphisms of the quaternions $\bH^r$ preserving a quaternionic hermitian form. It follows that they are Ricci flat. In fact, it follows from the theorems of De Rham and Berger, the Calabi-Yau theorem and results of Cheeger-Gromoll and Bochner that, after possibly taking a finite \'etale cover, Ricci-flat compact Riemannian manifolds are products of complex tori, Calabi-Yau manifolds and hyperk\"ahler manifolds (see Paragraph \ref{subsecdecomp}).

Constructing examples of compact hyperk\"ahler manifolds has proven particularly challenging. Two infinite series were constructed by Beauville, using an idea of Fujiki. Two sporadic families of hyperk\"ahlers of dimensions $6$ and $10$ were constructed by O'Grady (\cite{O'Grady1999}, \cite{O'Grady2003}) via desingularization of certain singular moduli spaces of sheaves on K3 surfaces and complex tori of dimension $2$. We give an overview of Beauville's constructions of the two infinite series.

It is the content of the Torelli theorem that hyperk\"ahler manifolds are essentially determined by their second cohomology. This is consistent with the fact that all constructions to date of hyperk\"ahler manifolds involve surfaces.

We briefly describe the moduli spaces of compact hyperk\"ahler manifolds, their period domains and some of their properties. By a result of Tian-Todorov and Bogomolov, the deformations of hyperk\"ahler manifolds are unobstructed. This essentially means that the moduli spaces of compact hyperk\"ahler manifolds are smooth analytic spaces. It is known however, that they are not Hausdorff.

The period domain of a given family of hyperk\"ahler manifolds is constructed from the lattice abstractly isometric to the second integral cohomology of the hyperk\"ahler together with a natural non-degenerate quadratic form called the Beauville-Bogomolov form. This form generalizes the intersection form in the case of dimension $2$ and the natural form on the second cohomology of the Fano variety of lines of a smooth cubic fourfold. In the case of the Fano variety of lines, the form is induced by the intersection form on the fourth cohomology of the cubic fourfold, via the Abel-Jacobi isomorphism between the second cohomology of the Fano variety if lines and the fourth cohomology of the cubic fourfold.

For a fixed compact hyperk\"ahler $X$, we describe the local and the global period domains with their respective maps from the local and global deformation spaces of $X$. We explain the local Torelli theorem and Verbitsky's weaker version of global Torelli which holds in the hyperk\"ahler case.

We conclude with a brief discussion of twistor conics and twistor families, the proof of the global Torelli theorem by Verbitsky and the relation between twistor families and hyperholomorphic bundles.

Some good general references for the material that we present here are: \cite{Beauville1983chern}, \cite{Beauville2007hol}, \cite{Beauville2011symp}, \cite{DoCarmo1992Riem}, \cite{GrossHuybrechtsJoyce2003}, \cite{VerbitskyKaledin1999}.

\section*{Acknowledgements}

This is an expanded version of the notes for a series of ten 45 minute lectures that I gave at the Trieste algebraic geometry summer school in July 2021. Exercise sessions for the lectures were run by Samir Canning, Yajnaseni Dutta and David Stapleton whom I wish to thank for their help. I also wish to thank the organizers: Valentina Beorchia, Ada Boralevi and Barbara Fantechi, for the invitation to lecture at the summer school and for the excellent organization of an enjoyable summer school, especially with the challenge of COVID-19.


\section{$C^\infty$ manifolds}\label{secinfty}

\subsection{Tangent and cotangent bundles}

For a $C^\infty$ manifold $M$, we denote by $T_M$ the tangent bundle of $M$ and $T_M^*$ the cotangent bundle.

For any non-negative integers $(k,l)$, the sections of the bundle $T_M^{\otimes k} \otimes (T_M^*)^{\otimes l}$ are called $(k,l)$-tensors. Section of $T_M$ are vector fields and sections of $\Lambda^p T_M^*$ differential $p$-forms. Alternatively, vector fields can be defined as first order differential operators on $C^\infty$ functions.

In a local coordinate chart with local coordinates $(x^1, \ldots, x^n)$, the (local) vector fields $\del /\del x^1, \ldots, \del /\del x^n$ form a basis of vector fields and the (local) $1$-forms $dx^1, \ldots, dx^n$ form a basis of differential $1$-forms. A local $(k,l)$-tensor can be written as
\[
T = \sum T^{i_1, \ldots, i_k}_{j_1, \ldots, j_l} \frac{\del}{\del x^{i_1}} \otimes \ldots \otimes \frac{\del}{\del x^{i_k}} \otimes dx^{j_1} \otimes \ldots \otimes dx^{j_l}.
\]

\subsection{The Lie bracket}

Given a vector field $v= \sum v^i \frac{\del}{\del x^{i}}$ and a $C^\infty$ function $f$ on $M$,
\[
v(f) = \sum_{i=1}^n v^i \frac{\del f}{\del x^i}.
\]
Given two vector fields $v= \sum v^i \frac{\del}{\del x^{i}}, w = \sum w^i \frac{\del}{\del x^{i}}$, the Lie bracket of $v$ and $w$ is given by
\[
[v,w] = \sum_{j=1}^n \left( \sum_{i=1}^n v^i \frac{\del w^j }{\del x^i} - w^i \frac{\del v^j }{\del x^i } \right) \frac{\del }{\del x^j }.
\]
Alternatively, the Lie bracket can be defined via its action on $C^\infty$ functions on $M$:
\[
[v,w] (f) = v (w(f)) - w (v(f)).
\]

\subsection{Connections}

Tangent vectors allow us to take derivatives of $C^\infty$ functions. Connections allow us to take derivatives of sections of arbitrary vector bundles.

 For a $C^\infty$ vector bundle $E$ on $M$, a connection is a linear map
\[
\nabla : C^\infty (E) \lra C^\infty (E \otimes T_M^*),
\]
satisfying the Leibnitz rule
\[
\nabla (fe) = f \nabla (e) + e \otimes df
\]
for all $C^\infty$ sections $e$ of $E$ and $C^\infty$ functions $f$ on $M$. For any vector field $v$ on $M$, the connection $\nabla$ defines a linear map $\nabla_v : C^\infty (E) \ra C^\infty (E)$ via
\[
\nabla_v (e) := \nabla (e)(v).
\]
We call $\nabla_v$ the covariant derivative in the direction of $v$.

We may thus also think of $\nabla$ as a linear map
\[
\nabla : C^\infty (E \otimes T_M) \lra C^\infty (E).
\]
When $E=T_M$, the torsion of a connection $\nabla : C^\infty (T_M \otimes T_M) \ra C^\infty (T_M)$ is the linear map
\[
T : C^\infty (\Lambda^2 T_M ) \lra C^\infty (T_M)
\]
defined as
\[
T (v \wedge w) := \nabla_v (w) - \nabla_w (v) - [v,w].
\]
We say $\nabla$ is torsion-free or symmetric when $T=0$.

\subsection{Curvature}

Euclidean space is ``flat''. What this means is that when we take second partial derivatives of vector fields, the order of differentiation does not affect the final result. Roughly speaking, the curvature of a connection measures the difference between the second partials of a section of a vector bundle taken in different orders. 

For general vector fields $v,w$, the curvature measures the difference between $\nabla_v\nabla_w - \nabla_w\nabla_v$ and the derivative in the direction of the bracket $[v,w]$. On the tangent bundle of Euclidean space this difference is $0$.

Precisely, the curvature of a connection $\nabla$ is, a linear map
\[
R : C^\infty (E) \lra C^\infty (E \otimes \Lambda^2T_M^*)
\]
or, equivalently, 
\[
R : C^\infty (E \otimes \Lambda^2 T_M) \lra C^\infty (E)
\]
or a global section
\[
R \in C^\infty (\End (E) \otimes \Lambda^2 T_M^*).
\]
It can be defined via its action on sections $e$ of $E$ and vector fields $v,w$ as
\[
R (e \otimes (v\wedge w)) = \nabla_v (\nabla_w(e)) - \nabla_w(\nabla_v(e)) - \nabla_{[v,w]} (e).
\]
We say that the connection $\nabla$ (or sometimes the bundle $E$) is flat if $R=0$.

In a coordinate chart with coordinates $(x^1, \ldots, x^n)$, the partial derivatives commute, i.e.,
\[
\left[\frac{\del}{\del x^i} , \frac{\del}{\del x^j }\right] = 0
\]
for all $i,j$. Hence
\[
R \left(e \otimes \left(\frac{\del}{\del x^i}\wedge \frac{\del}{\del x^j}\right)\right) = \nabla_\frac{\del}{\del x^i} \left(\nabla_\frac{\del}{\del x^j}(e)\right) - \nabla_\frac{\del}{\del x^j}\left(\nabla_\frac{\del}{\del x^i}(e)\right)
\]
and the connection is flat if and only if its partial (covariant) derivatives commute.

\subsection{Parallel transport}

Suppose given a $C^\infty$ vector bundle $E$ on $M$ with a connection
\[
\nabla : E \lra E \otimes T_M^*,
\]
and a smooth curve $\gamma : [0,1] \ra M$. Parallel transport along $\gamma$ produces sections of the pull-back $\gamma^* E$ that are `constant' or `horizontal' along $\gamma$. As we see below, such sections exist and are determined by their values at one point of $\gamma$. 

The pull-back $\gamma^* E$ is a $C^\infty$ vector bundle on $[0,1]$ with fiber $E_{\gamma (t)}$ at $t\in [0,1]$. The connection $\nabla$ defines the connection $\gamma^* \nabla$ on $\gamma^*E$ as the composition
\[
\gamma^*\nabla : \gamma^* E \lra \gamma^* E \otimes \gamma^* T_M^* \surj \gamma^* E \otimes T_{[0,1]}^*
\]
where the second map is induced by the projection $T_M^* \surj T_{[0,1]}^*$.

In local coordinates $(x^1, \ldots, x^n)$ on $M$, with $\gamma (t) = (x^1 (t), \ldots, x^n (t))$,
\[
\dot{\gamma} (t) = (\dot{x}^1 (t), \ldots, \dot{x}^n (t)) = \sum_{i=1}^n \dot{x}^i (t) \frac{\del}{\del x^i}
\]
and, for all (local) sections $e$ of $E$,
\[
\nabla_{\dot{\gamma} (t)} (e) := \nabla_{\sum_{i=1}^n \dot{x}^i (t) \frac{\del}{\del x^i}} (e) := \sum_{i=1}^n \dot{x}^i (t) \nabla_{\frac{\del}{\del x^i}} (e).
\]
\begin{defprop}
Put $x := \gamma (0), y := \gamma (1)$. Then, for all $e\in E_x = (\gamma^* E)_0$, there exists a unique smooth section $s$ of $\gamma^* E$ such that $s(0) =e$ and $\gamma^*\nabla (s) =0$, i.e., $\nabla_{\dot{\gamma}(t)} (s) =0$.

The parallel transport of $e$ along $\gamma$ to $y$ is $P_{\gamma} (e) := s(1) \in E_y = (\gamma^* E)_1$. The map
\[
P_{\gamma} : E_x \lra E_y
\]
is a linear isomorphism.
\end{defprop}

\subsection{Holonomy}

As we saw above, parallel transport defines linear isomorphisms between fibers of $E$ at points of $M$. In particular, for a given point $x$ of $M$, it defines linear automorphisms of the fiber $E_x$. The holonomy of $\nabla$ is the group generated by these automorphisms. It acts on all tensors of $E$ and its invariants are the covariantly constant tensors:

\begin{defprop}
If $\gamma$ is a loop (i.e. $x=y$), then $P_{\gamma} \in GL (E_x)$. The holonomy group $\Hol_x (\nabla)$ at $x$ is
\[
\Hol_x (\nabla ) := \{ P_\gamma \mid \gamma \hbox{ is a loop based at } x\}.
\]
It has the following properties.
\begin{enumerate}
\item $\Hol_x (\nabla )$ is a Lie subgroup of $GL (E_x)$:
\[
\gamma \delta (t) = \left\{ \begin{array}{lll}
\delta (2t) & \hbox{if} & t\in \left[ 0, \frac{1}{2} \right] \\
\gamma (2t -1) & \hbox{if} & t\in \left[ \frac{1}{2} , 1\right]
\end{array} \right.
\]
\[
\gamma^{-1} (t) = \gamma (1-t),
\]
\[
P_{\gamma \delta } = P_\gamma \circ P_\delta , \quad P_{\gamma^{-1}} = P_\gamma^{-1}.
\]
\item If $\gamma$ is a path from $x$ to $y$, then
\[
\Hol_y (\nabla) = P_\gamma \Hol_x (\nabla) P_\gamma^{-1}.
\]
Hence, up to conjugation, $\Hol_x (\nabla)$ only depends on the connected component of $M$ containing $x$.
\item if $M$ is simply connected, then $\Hol_x (\nabla)$ is connected. Any loop can be shrunk to a point:
\[
\gamma : [0,1] \times [0,1]  \lra M \: ; \quad \gamma_s (t) := \gamma (s, t) \: ; \gamma_1 (t) = x \hbox{ for all } t.
\]
Then $\{ P_s := P_{\gamma_s} \mid s\in [0,1] \}$ is a path in $\Hol_x (\nabla)$ from $P_0 = P_{\gamma_0}$ to $P_1 = P_{\gamma_1} = Id$.
\item Let $\Fhol_x (\nabla) \subset \Fgl (E_x) = \End (E_x)$ be the Lie algebra of $\Hol_x (\nabla)$. Recall that the curvature operator $R(\nabla)$ belongs to $C^\infty (E^* \otimes E \otimes \Lambda^2 T_M^* ) = C^\infty (\End (E) \otimes \Lambda^2 T_M^* )$. At a point $x$, the fiber $R(\nabla)_x$ of $R(\nabla)$ belongs to $\End (E_x) \otimes \Lambda^2 T_x^* M$. We have
\[
R(\nabla)_x \in \Fhol_x (\nabla) \otimes \Lambda^2 T_x^* M.
\]
\end{enumerate}
\end{defprop}

As we shall see below, Riemannian holonomy plays a central role in the structure theory of Riemannian manifolds.

The connection $\nabla$ induces connections on all tensor powers $E^{\otimes k} \otimes (E^*)^{\otimes l}$, and all exterior and symmetric powers of $E$ and $E^*$ and their tensor products. We shall denote these induced connections by $\nabla$ as well.
\begin{definition}
A tensor $S$ is called (covariantly) constant if $\nabla (S) =0$.
\end{definition}

\begin{theorem}
For a tensor $S$,
$\nabla (S) =0$ if and only if $S$ is fixed by $\Hol_x (\nabla)$, if and only if $P_\gamma (S(x)) = S(y)$ for all $x,y\in M$ and all paths $\gamma$ from $x$ to $y$.
\end{theorem}

\section{Riemannian manifolds}\label{secRiem}

A $C^\infty$ manifold is called Riemannian if it has a Riemannian metric, i.e., a (2,0)-tensor $g\in C^\infty ((T_M^*)^2$ which is symmetric:
\[
g\in C^\infty (\Sym^2 T_M^*),
\]
and defines a positive definite quadratic form on the tangent space $T_{M, x}$ for all $x\in M$. It is a fundamental result in differential geometry that every smooth manifold can be endowed with a Riemannian metric.

Riemannian manifolds have canonical connections on their tangent bundles: the Levi-Civita connection. The holonomy of the Levi-Civita connection is called Riemannian holonomy and the classification of Riemannian manifolds is based on the classification of Riemannian holonomy groups.

\subsection{Levi-Civita connection}

Suppose $(M,g)$ is a Riemannian manifold. The fundamental theorem of Riemannian geometry is the following.
\begin{theorem}
There exists a unique torsion free (or symmetric) connection $\nabla$ on $T_M$ such that $\nabla g=0$. This unique connection is called the Levi-Civita or Riemannian connection of $(M,g)$.
\end{theorem}

One can verify that the condition $\nabla g=0$ is equivalent to the following compatibility property: For all vector fields $u,v,w$ on $M$,
\[
u (g (v,w)) = g(\nabla_u v, w) + g (v, \nabla_u w).
\]
The Levi-Civita connection $\nabla$ can be explicitly defined via
\[
2 g (\nabla_u v, w) = u (g(v,w)) + v(g(u,w)) - w(g(u,v)) + g ([u,v],w) - g([v,w],u) - g([u,w], v).
\]
The curvature $R(\nabla)$ is a $(1,3)$ tensor:
\[
R(\nabla) : T_M \lra T_M \otimes \Lambda^2 T_M^*.
\]
More symmetries of $R(\nabla)$ can be exhibited by defining the $(0,4)$ tensor $\tR (\nabla)$ as the compostion
\[
\tR (\nabla ) : T_M \stackrel{R(\nabla)}{\lra} T_M \otimes \Lambda^2 T_M^* \stackrel{g \otimes Id}{\lra} T_M^* \otimes \Lambda^2 T_M^*.
\]
While a priori $\tR (\nabla) \in C^\infty ((T_M^*)^{\otimes 2} \otimes \Lambda^2 T_M^*)$, one can show that in fact
\[
\tR (\nabla )\in C^\infty (\Sym^2 ( \Lambda^2 T_M^*)).
\]
The Bianchi identities can be written in the form
\[
R(u,v) w + R(v,w) u+ R(w,u) v  =0, \quad \nabla_u R (u,v) + \nabla_v R (w,u) + \nabla_w R (u,v) =0.
\] 
In a basis of local coordinates $x^1, \ldots, x^n$, we can write $\tR (\nabla)$ as
\[ 
\tR (\nabla) = \sum_{a,b,c,d} \tR_{abcd} dx^a \wedge dx^b \odot dx^c \wedge dx^d,
\]
where $\alpha \odot \beta := \alpha \otimes \beta + \beta \otimes \alpha$ is the symmetric tensor.
The Bianchi identities then can be written as
\[
\tR_{abcd} + \tR_{acdb} + \tR_{adbc} = 0, \quad \frac{\del}{\del x^e} \tR_{abcd} + \frac{\del}{\del x^c} \tR_{abde} + \frac{\del}{\del x^d} \tR_{abec} =0.
\]

\subsection{Ricci curvature}

The Ricci curvature is a $(0,2)$ tensor, obtained by contracting $R(\nabla)$:

At each point $x\in M$, the curvature tensor $R$ defines a multilinear map
\[
\begin{array}{ccc}
R_x : T_x M \times T_x M \times T_x M & \lra & T_x M \\
(u,v,w) & \longmapsto & R(u,v) w
\end{array}
\]
The Ricci curvature is the $(0,2)$ tensor defined as
\[
\begin{array}{ccc}
\Ric_x : T_x M \times T_x M & \lra & \bR \\
(u,v) & \longmapsto & \trace (w \mapsto R_x (u,w)v)
\end{array}
\]
where $\trace$ is the trace of a linear map. It follows from the symmetries of the curvature tensor that the Ricci curvature is symmetric. In local coordinates, if we write the curvature tensor as
\[
R (\nabla) = \sum_{a,b,c,d} R^a_{bcd} \frac{\del}{\del x^a} \otimes dx^b \otimes dx^c \wedge dx^d,
\]
then the coordinates of the Ricci tensor are
\[
\Ric_{ab} = \sum_c R^c_{acb}.
\]
\begin{definition}
We say $g$ is an Einstein metric if the Ricci curvature is a constant multiple of the metric.
We say $g$ is Ricci flat if the Ricci curvature is $0$.
\end{definition}

\subsection{Riemannian holonomy}

For a Riemannian manifold $(M,g)$, the holonomy of the Levi-Civita connection $\nabla$ is called Riemannian holonomy. For $x\in M$, we write
\[
\Hol_x (g) := \Hol_x (\nabla)\subset \GL (T_x M), \quad \Fhol_x (g) := \Fhol_x (\nabla) \subset \Fgl (T_x M) = \End (T_x M) = T_x M \otimes T_x^* M.
\]
A first symmetry property of Riemannian holonomy is seen using the isomorphism
$g : T_M \ra T_M^*$.
\begin{proposition}
We have
\[
(g_x \otimes Id_x) (\Fhol_x (g) \subset \Lambda^2 T_x^* M.
\]
\end{proposition}
We saw that the curvature tensor $\tR \in (\Fhol_x (g) \otimes \Lambda^2 T_x^* M) \cap \Sym^2 (\Lambda^2 T_x^* M)$. Hence
\begin{theorem}
\[
\tR \in \Sym^2 \Fhol_x (g) \subset \Sym^2 (\Lambda^2 T_x^* M).
\]
\end{theorem}

\subsection{Reducibility}

The first step in the classification of Riemannian manifolds is to decompose them into their `irreducible' factors. As we see below, these correspond to the irreducible summands in the representation of the Riemannian holonomy group on the tangent space of $M$.

\begin{definition}
A Riemannian manifold is called (locally) reducible if every point has a neighborhood isometric to a product. It is called irreducible if it is not locally reducible. We have
\end{definition}
\begin{proposition}
Suppose a neighborhood of $x\in M$ is isometric to the product $(M_1, g_1) \times (M_2, g_2)$. Then
\[
\Hol_x (g_1 \times g_2) = \Hol_x (g_1) \times \Hol_x (g_2).
\]
\end{proposition}
\begin{theorem}
If $(M,g)$ is irreducible at $x$, then $\bR^n = T_x M$ is an irreducible representation of $\Hol_x (g)$.
\end{theorem}

\subsection{Symmetric and locally symmetric spaces}

A large and relatively well understood class of irreducible Riemannian manifolds is that of locally symmetric spaces.

\begin{definition}
A Riemannian manifold is called symmetric if, for all $p\in M$, there exists an isometry $s_p : M \ra M$ such that $s_p^2 = Id_M$ and $p$ is an isolated fixed point for $s_p$.
\end{definition}
\begin{definition}
A Riemannian manifold is called locally symmetric if every point has an open neighborhood isometric to an open subset of a symmetric space. It is called nonsymmetric if it is not locally symmetric.
\end{definition}
\begin{theorem}
$(M,g)$ is locally symmetric if and only if $\nabla R =0$.
\end{theorem}

\subsection{Geodesics and completeness}

To better understand locally symmetric spaces, we use `geodesics'. Geodesics allow us to define a notion of `completeness' (often called geodesic completeness) for Riemannian manifolds. Among other things, these notions allow us to describe all symmetric spaces in terms of Lie groups.

\begin{definition}
A geodesic is a parametrized smooth curve $\gamma : (a,b) \ra M$ such that, for all $t\in (a,b)$, $\nabla_{\dot{\gamma}(t)} \dot{\gamma} (t) =0$.
\end{definition}
Intuitively, a geodesic is the trajectory of a particle moving with constant velocity on the manifold: the equation $\nabla_{\dot{\gamma}(t)} \dot{\gamma} (t) =0$ means that the acceleration of the particle is $0$ with respect to the Levi-Civita connection.

The Riemannian metric defines a norm in the tangent space at each point of $M$. By integrating the length of the velocity vector of a parametrized (piecewise) smooth curve, we define the length of such a curve. One can show that geodesics are {\em locally} the `shortest' curves on $M$ for the Riemannian length. It can happen however that there are many geodesics of different lengths between two given points on a manifold. The simplest example of this is the cylinder with Riemannian metric induced from $\bR^3$. The Riemannian distance is defined as the infimum of the lengths of the (piecewise) smooth curves between two points on $M$. We have the following useful existence and uniqueness theorem for geodesics.
\begin{theorem}
For all $p\in M, v\in T_p M$, there exists a unique geodesic $\gamma : (a,b) \ra M$ such that $\gamma (0) =p, \dot{\gamma} (0) =v$.
\end{theorem}
\begin{definition}
A manifold $(M,g)$ is (geodesically) complete if every geodesic $(a,b) \ra M$ can be defined on all of $\bR \supset (a,b)$.
\end{definition}
All compact Riemannian manifolds and all symmetric spaces are complete. Every path connected Riemannian manifold which is also a complete metric space with respect to the Riemannian distance is geodesically complete.

We can now give the description of symmetric spaces in terms of Lie groups.
\begin{proposition}
Suppose $(M, g)$ is a connected, simply connected symmetric space. Then $(M, g)$ is complete.
Put
\[
G := \{ s_p \circ s_q \mid p,q\in M\} \subset \Isom (M).
\]
Then $G$ is a connected Lie group.
Choose $p\in M$ and let $H$ be the stabilizer subgroup of $p$ in $G$. Then $H$ is a closed connected Lie subgroup of $G$ and the map
\[
\begin{array}{ccc}
G/H & \lra & M \\
g & \longmapsto & g (p)
\end{array}
\]
is a diffeomorphism.
\end{proposition}

\subsection{De Rham's theorem}

De Rham's theorem describes the decomposition of a Riemannian manifold into the product of its irreducible factors.

\begin{theorem}
Suppose $(M,g)$ is Riemannian, complete, simply connected. Then $M$ is isometric to a product $M_0 \times M_1 \times \ldots \times M_k$ where $M_0$ is a Euclidean space, $M_1, \ldots, M_k$ are irreducible. The decomposition is unique up to reordering $M_1, \ldots, M_k$. The holonomy group of $M$ is the product of the holonomies of $M_1, \ldots, M_k$.
\end{theorem}

\subsection{Berger's theorem}

Suppose $(M, G)$ is connected. Then, $\Hol (g) := \Hol_x (g)$ is independent of the choice of $x$ up to conjugation in $\GL_n(\bR)$.
\begin{definition}
The restricted holonomy group $\Hol (g)^0$ is the connected component of the identity of $\Hol (g)$.
\end{definition}
 Berger's theorem classifies the possibilities for the restricted holonomy group $\Hol(g)^0$ and describes the corresponding manifolds.
\begin{theorem}
Suppose $(M,g)$ is Riemannian, complete, connected, nonsymmetric, irreducible. Then the restricted holonomy group $\Hol(g)^0$ is one of the following:
\begin{enumerate}
\item $\Hol(g)^0 \cong SO(n)$ (automorphisms of $\bR^n$ preserving the metric, generic metric),
\item $n = 2m \geq 4$, $\Hol(g)^0 = U(m) \subset SO(n)$ (automorphisms of $\bC^m$ perserving a hermitian form, K\"ahler),
\item $n = 2m \geq 4$, $\Hol(g)^0 = SU(m) \subset SO(n)$ (automorphisms of $\bC^m$, Calabi-Yau, Ricci-flat, K\"ahler),
\item $n = 4r \geq 4$, $\Hol(g)^0 = Sp(r) \subset SO(n)$ ($\bR$-linear automorphisms of $\bH^r$ preserving a quaternionic hermitian form, hyperk\"ahler, Ricci-flat, K\"ahler), (when $r=1$, the group $Sp(1)$ is abstractly isomorphic to the group $SU(2) = S^3$ of unit quaternions)
\item $n = 4r \geq 8$, $\Hol(g)^0 = Sp(r)Sp(1) \subset SO(n)$ ($\bR$-linear automorphisms of $\bH^r$, quaternionic-K\"ahler, Einstein, not Ricci-flat, not K\"ahler), (the group $Sp(1) = SU(2) = S^3$ of unit length quaternions acts on $\bH^r$ by right scalar multiplication and commutes with $Sp(r)$, however, this action is different from the action of $Sp(1)$ on $\bH$ preserving a quaternionic hermitian form; the Lie group $Sp(r)Sp(1)$ generated by combining this action with that of $Sp(r)$ is abstractly isomorphic to $(Sp(r) \times Sp(1) )/(\bZ /2 \bZ)$; when $r=1$, $Sp(1)Sp(1) = SO(4)$),
\item $n=7$, $\Hol(g)^0 = G_2 \subset SO(7)$ (automorphisms of $\Img \bO \cong \bR^7$, exceptional, Ricci-flat),
\item $n=8$, $\Hol(g)^0 = Spin (7) \subset SO(8)$ (automorphisms of $\bO \cong \bR^8$, exceptional, Ricci-flat).
\end{enumerate}
\end{theorem}

\section{K\"ahler manifolds}

For a complex manifold $M$, multiplication by $i$ defines an endomorphism $I : T_M \ra T_M$ satisfying $I^2 =-Id$. This is called the complex structure (operator) of $M$. A metric $g$ on $M$ is called Hermitian if
\[
g(v,w) = g(Iv, Iw), \quad \hbox{for all vector fields} \quad v,w.
\]
The $(1,1)$ form associated to $g$ and $I$ is
\[
\omega (v,w) := g (Iv, w), \quad \hbox{for all vector fields} \quad v,w.
\]
Equivalently, $\omega$ is the composition
\[
\omega : T_M \stackrel{I}{\lra} T_M \stackrel{g}{\lra} T_M^*.
\]
The fact that $\omega$ is a $(1,1)$ form means $\omega (Iv, Iw) = \omega(v,w)$. One also checks that $\omega$ is anti-symmetric.

It is easy to check that any two of $\{I, g, \omega \}$ determine the third.
\begin{defprop}
The metric $g$ is K\"ahler with respect to $I$ if one of the following equivalent conditions hold:
\begin{enumerate}
\item $d \omega =0$,
\item $\nabla \omega =0$,
\item $\nabla I =0$.
\end{enumerate}
In such a case, $\omega$ is called the K\"ahler form of $g$.
\end{defprop}
So $g$ is K\"ahler if and only if $\omega$ and $I$ are constant. Equivalently $\Hol(g)$ preserves $\omega$ and $I$. The subgroup of $SO(n)$ preserving $I$ is $U(m)$ ($n=2m$). Therefore, $M$ is K\"ahler if and only if $\Hol (g) \subset U(m)$.

\subsection{Ricci form}

Given a K\"ahler manifold $(M,g,I)$, its Ricci form $\rho$ is the differential form associated to the Ricci curvature via $I$:
\[
\rho (v,w) := \Ric (Iv, w), \quad \hbox{for all vector fields} \quad v,w.
\]
Equivalently, $\rho$ is the composition
\[
\rho : T_M \stackrel{I}{\lra} T_M \stackrel{\Ric}{\lra} T_M^*.
\]
As in the case of $\omega$ and $g$: $\rho \in C^\infty (\Lambda^2 T_M^*)$. We have the following
\begin{proposition}
The Ricci form $\rho$ is a closed $(1,1)$ form. Its cohomology class in $H^2 (M, \bR)$ is $[\rho] = 2\pi c_1 (K_M) = 2\pi c_1 (T_M^*)$.
\end{proposition}

\subsection{Ricci flatness (the Calabi-Yau case)}

The Ricci form is the curvature of the connection induced on $K_M := \Omega^m_M$ by the Levi-Civita connection. So, if $\rho =0$, then $K_M$ is a flat bundle.

Assume now that $M$ is Ricci-flat and simply connected. The flat bundle $K_M$ admits local flat, i.e., covariantly constant, sections. Since $M$ is simply connected, $K_M$ has a global flat section. Such a section is hence invariant under Riemannian holonomy and, by the following lemma which is a consequence of Bochner's principle, holomorphic.
\begin{lemma}
Suppose $(M,I,g)$ is a compact K\"ahler, simply connected, Ricci-flat manifold with holonomy group $H$. For all $x\in M$ and all positive integers $p$, the natural evaluation map
\[
\begin{array}{ccc}
H^0 (M, \Omega^p_M) & \lra & (\Omega^p_{M,x})^H \\
w & \longmapsto &w_x
\end{array}
\]
is an isomorphism.
\end{lemma}
Hence $K_M$ has a nowhere vanishing holomorphic section, which implies that $K_M$ is trivial, i.e., $M$ is Calabi-Yau. Furthermore, on the tangent space $T_p M$ at a point $p\in M$, a nonvanishing differential $m$-form is a multiple of the determinant. Hence $\Hol (g)$ preserves the determinant. Since we already know that $\Hol (g) \subset U(m)$, this implies that $\Hol (g) \subset SU(m)$.

Conversely, if $\Hol (g) \subset SU(m)$, then $M$ admits a nowhere vanishing differential $m$-form, $K_M$ is trivial and $\rho =0$.

\subsection{The hyperk\"ahler case}

Recall that the quaternions have bases of the form
\[
\bH = \bR 1 \oplus \bR i \oplus \bR j \oplus \bR k, \quad \hbox{with} \quad i^2=j^2=k^2=ijk=-1.
\]
A triple $(i,j,k)$ as above is called a quaternionic triple. The Lie group $Sp(r)$ is the group of $\bR$-linear endomorphisms of $\bH^r$ preserving a quaternionic Hermitian form $q$. Recall that $q$ is quaternionic Hermitian if
\[
q (av, bw) = \oa\, b\, q(v,w), \quad \hbox{for all } a,b\in \bH, v,w\in \bH^r
\]
where, if $a = \lambda + \mu i + \nu j + \rho k$, then $\oa = \lambda - \mu i - \nu j - \rho k$. Such a $q$ can be represented by an $r\times r$ matrix $A$ with entries in $\bH$ such that $A\, \oA^t = Id$ is the identity of $\bH^r$.

We can embed $Sp(r)$ in $SU(2r)$ each time we choose $i\in\bH$ with $i^2=-1$ as follows.

Complete $i$ to a quaternionic triple $(i,j,k)$ and write
\[
q = h + \omega j
\]
where $h$ is Hermitian with respect to $i$ and $\omega$ is alternating $\bC$-bilinear with respect to the complex structure on $\bH^r$ given by $i$. Then $Sp(r)$ can be identified with the group of $\bR$-linear automorphisms of $\bH$ preserving $h$ and $\omega$. Hence, thinking of $U(2r)$ as the group of transformations of $\bH^r = \bC \oplus \bC i$ preserving $h$, we can identify $Sp(r)$ as the subgroup of $U(2r)$ of transformations preserving $\omega$. In particular, they preserve $\wedge^r\omega$, which means they belong to $SU(2r)$.

Given a Riemannian manifold $M$ with $\Hol_p (g) \subset Sp(r)$, we can identify $T_p M$ with $\bH^r$. The form $\omega$ obtained as above by decomposing the form $q$ is invariant under the holonomy group of $M$, hence globalizes to an alternating flat, i.e., holomorphic, $2$-form on $M$ which is non-degenerate everywhere. Furthermore, the quaternionic triple $(i,j,k)$ gives three complex structures $I,J,K$ on $M$ satisfying the quaternionic relations and with respect to which $g$ is K\"ahler ($I,J,K$ are invariant under $\Hol_p (g)$, hence flat). We then obtain a sphere of complex structures $\lambda = aI + bJ + cK$ with $a,b,c\in \bR, a^2 + b^2 + c^2 =1$ such that $\nabla \lambda =0$. The metric $g$ is therefore K\"ahler with respect to all these complex structures.

Note that if $\Hol (g) = U(m)$ or $SU(m)$, then $M$ has a unique complex structure with respect to which $g$ is K\"ahler because the only complex endomorphisms commuting with $U(m)$ or $SU(m)$ are multiplication by scalars. So Calabi-Yaus have only one K\"ahler complex structure.

If $\Hol (g) = Sp(r)$, then $M$ has exactly an $S^2$ of K\"ahler complex structures because the only quaternionic endomorphisms commuting with $Sp(r)$ are multiplication by quaternionic scalars.

If $M$ is a complex torus, then $\Hol (g) =0$. Any complex structure is then K\"ahler.

\begin{definition}
We say that $M$ is irreducible hyperk\"ahler if $\Hol(g) = Sp(r)$, i.e., $M$ has exactly an $S^2$ of K\"ahler complex strcutures.
\end{definition}

\subsection{The Calabi conjecture and its consequence}

\begin{theorem} Calabi's conjecture, Yau's theorem:

Let $(M, I)$ be a compact complex manifold and $g$ a metric K\"ahler with respect to $I$ with K\"ahler form $\omega$ and Ricci form $\rho$. Let $\rho'$ be a real closed $(1,1)$ form on $M$ with cohomology class $[\rho'] = [\rho] = 2\pi c_1 (K_M)$. There exists a unique K\"ahler metric $g'$ on $M$ whose Ricci form is $\rho'$ and whose K\"ahler form $\omega'$ satisfies $[\omega'] = [\omega]$.
\end{theorem}
For Ricci-flat manifolds this has the following useful consequence.
\begin{corollary}
Suppose $(M, I, g)$ is compact K\"ahler with $c_1 (K_M) =0$. There exists a unique Ricci-flat K\"ahler metric in each K\"ahler class on $M$. The Ricci-flat K\"ahler metrics on $M$ form a smooth family of dimension $h^{1,1} (M)$, isomorphic to the K\"ahler cone of $M$.
\end{corollary}

\subsection{The decomposition theorem}\label{subsecdecomp}

The following decomposition theorem for Ricci-flat manifolds is a consequence of De Rham's decomposition theorem, the Berger classification theorem and results of Cheeger-Gromoll and Bochner. (see \cite[Th\'eor\`eme 1]{Beauville1983chern}).
\begin{theorem}
Let $(M, I, g)$ be a compact K\"ahler, complete, Ricci-flat manifold. Then
\begin{enumerate}
\item the universal cover of $M$ is isomorphic to $\bC^k \times \prod_i V_i \times \prod_j X_j$ where $\bC^k$ has the standard K\"ahler metric, and, for all $i$, $V_i$ is compact simply connected with holonomy $SU(m_i)$ and, for all $j$, $X_j$ is compact simply connected with holonomy $Sp (r_j)$,
\item there exists a finite \'etale cover of $M$ isomorphic to $T \times \prod_i V_i \times \prod_j X_j$ where $T$ is a complex torus of complex dimension $k$.
\end{enumerate}
\end{theorem}
The proof uses
\begin{lemma}
Suppose $(M,I,g)$ is a compact K\"ahler, simply connected, Ricci-flat manifold. The group of automorphisms of $(M,I)$ is discrete. In particular, the group of automorphisms of $(M,I,g)$ is finite (because it is contained in $SO(n)$ which is compact).
\end{lemma}

\section{Holomorphic symplectic manifolds}

We now present the infinite series of examples of compact hyperk\"ahler manifolds constructed by Beauville \cite{Beauville1983chern}. For this, the point of view of holomorphic symplectic geometry is more convenient. We begin with the following.

\begin{proposition}
Suppose $(M,I,g)$ is a compact K\"ahler, simply connected, Ricci-flat manifold of complex dimension $2r$ with holonomy group $Sp(r)$. Then
\begin{enumerate}
\item there exists a holomorphic $2$-form $\varphi$ on $M$ which is nondegenerate everywhere (represented by the form $\omega$ in the decomposition of the quaternionic Hermitian form $q = h + \omega j$),
\item for all $0\leq p\leq r$,
\[
H^0 (M, \Omega_M^{2p+1}) =0, \quad H^0 (M, \Omega_M^{2p}) = \bC \varphi^p.
\]
\end{enumerate}
\end{proposition}
\begin{defprop}
A compact K\"ahler manifold $X$ is called holomorphic symplectic if there exists an everywhere non-degenerate holomorphic $2$-form on $X$. This is equivalent to: $X$ is compact hyperk\"ahler or $X$ is K\"ahler and $Hol_g(X) \subset Sp(r)$.

A compact K\"ahler manifold $X$ is called irreducible holomorphic symplectic if $X$ is simply connected and $H^2 (X, \Omega^2_X)$ is generated by an everywhere non-degenerate holomorphic $2$-form. This is equivalent to: $X$ is irreducible compact hyperk\"ahler $X$ is K\"ahler and $Hol_g(X) = Sp(r)$.
\end{defprop}

\subsection{The case of surfaces}

In dimension $2$, $Sp(1) = SU(2)$, so Calabi-Yau and hyperk\"ahler are the same: these are K3 surfaces and complex tori.
\begin{definition}
A K3 surface is a compact complex manifold of dimension $2$ such that $\Omega_X^2 \cong \cO_X$ and $H^1 (X, \cO_X) =0$.
\end{definition}
One can prove that K3 surfaces are simply connected and their integral cohomology is torsion free.

It is a deep theorem of Siu that a K3 surface admits a unique K\"ahler metric.

Examples of algebraic K3 surfaces:
\begin{enumerate}
\item Double covers of $\bP^2$ branched along smooth sextics.
\item Smooth quartics in $\bP^3$.
\item $(2,3)$ complete intersections in $\bP^4$.
\item $(2,2,2)$ complete intersections in $\bP^5$.
\end{enumerate}

\subsection{Hilbert schemes of points}

Both infinite series of examples are constructed using the Hilbert schemes of points, the first uses the Hilbert schemes of points of K3 surfaces, and the second uses the Hilbert schemes of points of complex tori of dimension $2$. The construction begins by showing that these Hilbert schemes have natural holomorphic symplectic structures.

Suppose $S$ is a compact complex manifold of dimension $2$. Denote $S^r$ the $r$-th Cartesian power of $S$ and
\[
\pi : S^r \surj S^{(r)} := S^r /\FS_r
\]
its quotient by the action of $\FS$ permuting the factors. Let $\Delta_{ij} \subset S^r$ be the diagonal where the $i$-th and $j$-th components are equal. The action of $\FS_r$ is not free on the diagonals $\Delta_{ij}$. The stabilizer of a generic point of $\Delta_{ij}$ is the subgroup $\{ 1, (ij)\}\subset \FS_r$ where $(ij)$ is the transposition exchanging $i$ and $j$. The quotient morphism $\pi$ is \'etale away from $\cup_{i, j} \Delta_{ij}$. Since the diagonals $\Delta_{ij}$ have codimension $2$ in $S^r$, by the theorem on the purity of the ramification locus of a morphism of smooth varieties, the symmetric power $S^{(r)}$ is singular along the diagonal $D := \pi (\Delta_{ij}) = \pi (\cup_{i,j} \Delta_{ij})$. Note that $D$ is irreducible.

The symmetric power $S^{(r)}$ has a natural desingularization: the Hilbert scheme $S^{[r]}$ of length $r$ Artinian subschemes of $S$. The natural map $\epsilon : S^{[r]} \ra S^{(r)}$ sends a subscheme $Z$ of length $r$ to its underlying $0$-cycle. Since, for any $r$ distinct points $x_1, \ldots, x_r\in S$, there exists a unique Artinian subscheme supported on $\{x_1, \ldots, x_r\}$, the map $\epsilon : S^{[r]} \setminus \epsilon^{-1} (D) \ra S^{(r)} \setminus D$ is an isomorphism.

Let $D_*\subset D$ be the open subset where exactly two coordinate are equal. Given $2x_1 + x_2 + \ldots + x_{r-1} \in D_*$, the datum of an Artinian subscheme of length $r$ supported on $2x_1 + x_2 + \ldots + x_{r-1}$ is equivalent to the datum of a tangent line to $S$ at $x_1$. So the set of Artinian subschemes of length $r$ supported on $2x_1 + x_2 + \ldots + x_{r-1}$ is naturally identified with $\bP T_{x_1} S$.

Let $S^{(r)}_*$, respectively $S^r_*$, be the open subset where at most two of the coordinates coincide and let $S^{[r]}_*$ be the inverse image of $S^{(r)}_*$ in $S^{[r]}$. The fiber of $\epsilon : S^{[r]}_* \ra S^{(r)}_*$ at $x\in D_*$ is naturally identified with $\bP T_{x_1} S$. One can prove:
\begin{theorem}
\begin{enumerate}
\item The complex analytic pair $(S^{(r)}_*, D_*)$ is locally isomorphic to $(B \times C, B\times \{O\})$, where $B$ is a ball, $C$ is a cone with vertex $O$ over a smooth conic in $\bP^2$.
\item The complex manifold $S^{[r]}_*$ is the blow up of $S^{(r)}_*$ along $D_*$.
\item If we denote $Bl_\Delta (S^r_*)$ the blow up of $S^r_*$ along the union of its diagonals, then the action of $\FS_r$ lifts to $Bl_\Delta (S^r_*)$ and
\[
S^{[r]}_* = Bl_\Delta (S^r_*) /\FS_r.
\]
\end{enumerate}
\end{theorem}
So we have the Cartesian diagram
\[
\xymatrix{
         Bl_\Delta (S^r_*) \ar[d]_{\rho} \ar[r]^{\eta} & S^r_* \ar[d]_\pi \\
         S^{[r]}_* \ar[r]^{\epsilon}      & S^{(r)}_*.}
\]
Next we construct differential forms on $S^{[r]}$, starting from differential forms on $S$.

Given a holomorphic differential form $\omega$ on $S$, the form $\psi := pr_1^* \omega + \ldots + \pr_r^* \omega$ and its pull-back $\eta^* \psi$ to $Bl_\Delta (S^r_*)$ are invariant under the action of $\FS_r$. Hence there exists a holomorphic differential form $\varphi$ on $S^{[r]}_*$ such that
\[
\eta^* \psi = \rho^* \varphi.
\]
\begin{proposition}
If $K_S = \Omega^2_S$ is trivial, then $S^{[r]}$ admits a holomorphic symplectic form.
\end{proposition}
\begin{proof}
Let $\omega$ be a generator of $K_S$. Defining $\psi$ and $\varphi$ as above, we show that $\varphi$ extends to $S^{[r]}$ as an everywhere non-degenerate form.

The form $\varphi$ extends to all of $S^{[r]}$ because $S^{[r]} \setminus S^{[r]}_*$ has codimension $\geq 2$ in $S^{[r]}$. The fact that $\varphi$ is everywhere non-degenerate means that $\wedge^r \varphi$ does not vanish anywhere.

The form $\wedge^r \varphi$ is a section of $K_{S^{[r]}}$, so the locus where it vanishes is a canonical divisor on $S^{[r]}$.

Denote $E_{ij} := \eta^* \Delta_{ij}$. Then the divisors $E_{ij}$ are the exceptional divisors of the blow up $\eta : Bl_\Delta (S^r_*) \ra S^r_*$ and the ramification divisors of the morphism $\rho : Bl_\Delta (S^r_*) \ra S^{[r]}_*$. Hence
\[
K_{Bl_\Delta (S^r_*)} = \rho^* K_{S^{[r]}_*} + \sum_{i <j} E_{ij},
\]
and the divisor of zeros of $\rho^* \wedge^r \varphi$ is
\[
\Div (\rho^* \wedge^r \varphi) = \rho^* \Div (\wedge^r \varphi) + \sum_{i <j} E_{ij}.
\]
However,
\[
\Div (\rho^* \wedge^r \varphi) = \Div (\eta^* \wedge^r \psi) = \Div (\wedge^r \eta^* \psi) = \sum_{i <j} E_{ij}.
\]
Indeed, choose $z = (x_1, \ldots , x_r ) \in S^r$, then
\[
T_z S^r = T_{x_1} S \oplus \ldots \oplus T_{x_r} S.
\]
The differential form $\psi$ is a bilinear form on $T_z S^r$, the decomposition $T_z S^r = T_{x_1} S \oplus \ldots \oplus T_{x_r} S$ is orthogonal with respect to $\psi$ and $\psi$ is non-degenerate at any $z$. Hence $\Div (\wedge^r \psi) =0$ on $S^r$. However, the differential of the blow up $\eta : Bl_\Delta (S^r_*) \ra S^r_*$ has image of dimension $2r-1$ along the union of the diagonals, so $\eta^* \psi$ is degenerate of rank $2r-2$ along $\cup_{i<j} E_{ij}$. It follows that $\Div (\wedge^r \eta^* \psi) = \sum_{i<j} E_{ij}$.

So $\rho^* \Div (\wedge^r \varphi) =0$ and $\Div (\wedge^r \varphi) =0$.
\end{proof}

To determine the type of $S^{[r]}$, we compute its fundamental group. The map $S^r \ra S^{(r)}$ is a Galois cover with Galois group $\FS_r$. So we have the exact sequence of fundamental groups
\[
1 \lra \FS_r \lra \pi_1 (S^r) \lra \pi_1 (S^{(r)}) \lra 1.
\]
We have
\[
\pi_1 (S^r) = \pi_1 (S^r_*) = \pi_1 (BL_\Delta (S^r_*)), \quad \pi_1 (S^{(r)}) = \pi_1 (S^{(r)}_*), \quad \pi_1 (S^{[r]}) = \pi_1 (S^{[r]}_*).
\]
The map $BL_\Delta (S^r_*) \ra S^{[r]}_*$ is also a Galois cover with Galois group $\FS_r$. So we have the commutative diagram of exact sequences
\[
\xymatrix{
         1 \ar[r] & \FS_r \ar[r]^{} \ar@{^{}=}[d]_{} & \pi_1 (BL_\Delta (S^r_*)) \ar@{^{}=}[d]_{} \ar[r]^{} & \pi_1 (S^{[r]}_*) \ar[d] \ar[r] & 1 \\
          1 \ar[r] & \FS_r \ar[r] & \pi_1 (S^r) \ar[r] & \pi_1 (S^{(r)}) \ar[r] & 1. }
\]
Therefore, we also have $\pi_1 (S^{[r]}_*) \stackrel{\cong}{\ra} \pi_1 (S^{(r)})$.

It is a fact from algebraic topology and group theory that $\pi_1 (S^{(r)})$ is the largest commutative quotient of $\pi_1 (S)$, hence it is isomorphic to $H_1 (S, \bZ)$.
\begin{lemma}
\begin{enumerate}
\item \[ H^i (S^{(r)}, \bQ) = H^i (S^r, \bQ )^{\FS_r} \]
\item \[ H^2 (S^{[r]}, \bQ ) = H^2 (S^{(r)}, \bQ) \oplus \bQ [E] \]
\item \[ H^2 (S^{(r)}, \bQ) = H^2 (S, \bQ ) \oplus \Lambda^2 H^1 (S, \bQ) \]
\end{enumerate}
\end{lemma}
\begin{proof}
\begin{enumerate}
\item Standard.
\item Replace $S^r$ by $S^r_*$, $S^{(r)}$ by $S^{(r)}_*$ and $S^{[r]}$ by $S^{[r]}_*$: the second cohomology does not change. We compute
\[
H^2 (S^{[r]}_*, \bQ) = H^2 (BL_\Delta (S^r_*), \bQ)^{\FS_r} = \left(H^2 (S^r_*, \bQ) \oplus (\oplus_{1\leq i <j \leq r} \bQ [E_{ij}])\right)^{\FS_r} = H^2 (S^r_*,  \bQ)^{\FS_r} \oplus \bQ [\rho^* E].
\]
\item We compute, using part (1),
\[
\begin{split}
H^2 (S^{(r)}, \bQ ) = H^2 (S^r, \bQ )^{\FS_r} \cong \left( H^2 (S, \bQ )^{\oplus r} \oplus \left(H^1 (S, \bQ)^{\otimes 2}\right)^{\oplus {r \choose 2}}\right)^{\FS_r} \\
= H^2 (S, \bQ) \oplus \left(H^1 (S, \bQ)^{\otimes 2}\right)^{transposition} \cong H^2 (S, \bQ) \oplus \Lambda^2 H^1 (S, \bQ)
\end{split}
\]
by skew-symmetry.
\end{enumerate}
\end{proof}
We immediately obtain.
\begin{corollary}
If $S$ is a K3 surface, then $S^{[r]}$ is an irreducible holomorphic symplectic manifold and
\[
H^2 (S^{[r]}, \bQ ) = H^2 (S, \bQ) \oplus \bQ [E].
\]
\end{corollary}
$S^{[r]}$ is K\"ahler by results of Varouchas.

\subsection{Generalized Kummers}

Now take $S=A$ a complex torus of dimension $2$. Then $A^{[r+1]}$ is a holomorphic symplectic manifold. As in the case of K3 surfaces, it is K\"ahler. By the previous results,
\[
\pi_1 (A^{[r+1]} )= H_1 (A, \bZ) = \pi_1 (A) \neq \{1\}, \quad H^2 (A^{[r+1]}, \bQ) = H^2 (A, \bQ) \oplus \Lambda^2 H^1 (A, \bQ ) \oplus \bQ [E].
\]
So in this case, the Hilbert scheme is not irreducible holomorphic symplectic. We determine its factors according to the decomposition theorem.

Consider the addition map $s : A^{(r+1)} \ra A$ and its composition
\[
S : A^{[r+1]} \stackrel{\rho}{\lra} A^{(r+1)} \stackrel{s}{\lra} A.
\]
\begin{definition}
The $(r+1)$-st generalized Kummer manifold of $A$ is
\[
K_r := S^{-1} (0).
\]
\end{definition}
One can see that $K_r$ is a manifold as follows.

The complex torus $A$ acts on itself by translation, hence also on $A^{[r+1]}$ by pull-back:

If $Z\subset A$ is an analytic subspace of length $r+1$, then $a\in A$ acts as $Z \mapsto t_a^* Z$ on $A^{[r+1]}$. The map $S$ is equivariant for this action on $A^{[r+1]}$ and the action of $A$ on itself via $x \mapsto t_{(r+1)a}^*x$. In other words we have the Cartesian diagram
\[
\xymatrix{
         A \times A^{[r+1]} \ar[d]_{} \ar[r]_{(a,Z)\mapsto t_a^* Z} & A^{[r+1]} \ar[d]^{S} \\
         A \times A \ar[r]_{(a,x)\mapsto t_{(r+1)a}^* x}      & A}
\]
which induces the Cartesian diagram
\[
\xymatrix{
         A \times K_r \ar[d]_{} \ar[r]^{(a,Z)\mapsto t_a^* Z} & A^{[r+1]} \ar[d]^{S} \\
         A \ar[r]_{a\mapsto (r+1)a}      & A.}
\]
It follows that $S$ is a smooth map and all its fibers are isomorphic to $K_r$ which is therefore also smooth.
\begin{proposition}
The holomorphic symplectic structure of $A^{[r+1]}$ restricts to a holomorphic symplectic structure on $K_r$.
\end{proposition}
\begin{proof}
Since $K_r$ is a fiber of a smooth morphism, its normal bundle is trivial: the normal space at every point of $K_r$ maps isomorphically onto $T_0A$, so that we have $N_{K_r/ A^{[r+1]}} \cong T_0 A \otimes \cO_{K_r}$. From the normal bundle sequence
\[
0 \lra T_{K_r} \lra T_{A^{[r+1]}} |_{K_r} \lra N_{K_r/ A^{[r+1]}} \lra 0
\]
we obtain $K_{K_r} \cong K_{A^{[r+1]}} |_{K_r} \cong \cO_{K_r}$.

Recall the differential forms $\psi = pr_1^* \omega \oplus \ldots \oplus pr_{r+1}^* \omega$ and $\varphi$ with $\eta^* \psi = \rho^* \varphi$. The form $\wedge^r (\varphi |_{K_r})$ is a section of $K_{K_r} \cong \cO_{K_r}$. We show that it remains everywhere non-degenerate. As before, this means that $\wedge^r (\varphi |_{K_r})$ does not vanish anywhere. Since $K_{K_r}$ is trivial, either $\wedge^r (\varphi |_{K_r})$ is zero everywhere or it does not vanish anywhere. We prove that it is nonzero at one point.

Let $Z = x_1 + \ldots + x_{r+1}\in K_r$ be such that the $x_i$ are all distinct. Then
\[
T_Z A^{[r+1]} \cong T_(x_1, \ldots, x_{r+1}) A^{r+1} \cong T_{x_1} A \oplus \ldots \oplus T_{x_{r+1}} A \cong (T_0 A)^{\oplus (r+1)}.
\]
We can choose the isomorphism above in such a way that the differential $dS : T_Z A^{[r+1]} \ra T_0 A$ of $S$ is the sum map. The form $\varphi$ acts as $\omega$ on each summand $T_0A$ of $T_Z A^{[r+1]}$ and the summands are orthogonal to each for $\varphi$. It is then an exercise in linear algebra to check that $\varphi |_{\Ker dS}$ is non-degenerate, i.e., $\wedge^r (\varphi |_{K_r})$ is not $0$.
\end{proof}
\begin{proposition}
The manifold $K_r$ is simply connected. For $r\geq 2$, we have
\[
H^2 (K_r, \bQ) \cong H^2 (A, \bQ) \oplus \bQ [E]
\]
where $E$ is the intersection of the exceptional divisor of $A^{[r+1]}$ with $K_r$.
\end{proposition}
\begin{proof}
Immediate from the definition of $K_r$ and the description of the cohomology and fundamental group of $A^{[r+1]}$.
\end{proof}
It now follows that the factors of $A^{[r+1]}$ in the decomposition theorem are $K_r$ and $A$ itself.

Note that $S^{[r]}$ (for K3 surfaces $S$) and $K_r$ have different betti numbers, hence are not deformation equivalent. These provide two infinite series of families of hyperk\"ahler manifolds.

There are two known examples of families of hyperk\"ahler manifolds due to O'Grady that are not deformation equivalent to Hilbert schemes of K3s or generalized Kummers: these are hyperk\"ahlers of dimensions $6$ and $10$.

\begin{question}
Are there other families of compact irreducible hyperk\"ahlers?
\end{question}

\section{Moduli of hyperk\"ahlers, the Beauville-Bogomolov form, the period domain and the period map}\label{secApp}

\subsection{Moduli of complex structures and Teichm\"uller space}

Given a differentiable manifold $X$, there can be many different complex structures on $X$. We define the Teichm\"uller space of $X$ as
\[
\Teich (X) := \{ \hbox{complex structures on } X \} /\sim^0
\]
where two complex structures $I, J$ on $X$ satisfy $I\sim^0 J$ if there exists a diffeomorphism $\varphi :X \ra X$ isotopic (or homotopic) to the identity $Id_X$ such that $\varphi^*I = J$. The moduli space of complex structures on $X$ is, by definition,
\[
\cM_{cx} (X) := \{ \hbox{complex structures on } X \} /\sim
\]
where two complex structures $I, J$ on $X$ satisfy $I\sim J$ if there exists a diffeomorphism $\varphi :X \ra X$ such that $\varphi^*I = J$. If we denote $\Diff(X)$ the group of diffeomorphisms of $X$ and $\Diff^0 (X)$ its connected component of the identity, then $G := \Diff (X) / \Diff^0 (X)$ is the discrete group of components of $\Diff (X)$, and
\[
\cM_{cx} (X) = \Teich (X) / G.
\]
A priori, $\cM_{cx} (X)$ is the space that we are interested in. However, it usually does not have many good properties while $\Teich (X)$ does. So we will, most of the time, work with small open sets of $\Teich (X)$ which describe small deformations of given complex structures.

\subsection{Universal families and Kuranishi's theorem}

Suppose given a complex manifold $(X, I)$.
\begin{definition}
A family of complex manifolds is a smooth proper morphism of complex spaces
\[
\pi : \cX \ra S.
\]
A deformation of $(X,I)$ is a family of complex manifolds
with a point $s_0 \in S$ and an isomorphism $\cX_0 := \pi^{-1} (s_0) \cong X$.

A deformation is called universal if, for any deformation $\cX' \ra S'$, there exists a unique morphism $\varphi : S' \ra S$ such that $\varphi (s'_0) = s_0$ and $\cX' \ra S'$ is the pull-back of $\cX \ra S$ under $\varphi$. In other words, we have the Cartesian diagram
\[
\xymatrix{
         \cX' \ar[d]_{} \ar[r]^{} & \cX \ar[d]^{\pi} \\
         S' \ar[r]_{\varphi}      & S.}
\]
The universal deformation is unique up to unique isomorphism and we denote it
\[
\cX \ra \Def (X).
\]
\end{definition}
Kuranishi's theorem is the following.
\begin{theorem}
Suppose $(X,I)$ is a compaxt complex manifold with $H^0 (X, T_X) =0$. Then a \emph{local} universal deformation of $(X,I)$ exists and it is universal for all of its fibers.
\end{theorem}
Under the conditions of the theorem, the local universal deformation $\cX \ra \Def (X)$ is sometimes called the Kuranishi family.

Note that the condition $H^0 (X, T_X) =0$ means there are no global holomorphic vector fields on $X$ or $X$ has no infinitesimal automorphisms: given two complex manifolds $X, Y$ and a holomorphic map $f: X \ra Y$, the tangent space to the space of holomorphic maps $\Hom (X,Y)$ at $f$ can be identified with $H^0 (X, f^* T_Y)$. This can be deduced from general results in deformation theory, applied to the deformations of the graph of $f$ in $X\times Y$.

\subsection{Unobstructedness for $K$-trivial K\"ahler manifolds}

For any compact complex manifold $X$, if $H^0 (X, T_X) =0$, then $X$ has a local or small universal deformation denoted $\cX \ra \Def (X)$. By this we mean a germ of a deformation, i.e., whose base is suitably small. Such a deformation is universal for all its fibers, its base $\Def (X)$ is a ``Kuranishi slice'' $\subset H^1 (X,T_X)$. For $t\in \Def (X)$ small, we have
\[
T_t \Def (X) = H^1 (X_t, T_{X_t}).
\]
The obstructions to deformations (to various orders) provide local analytic equations for $\Def (X)$ in a neighborhood of $0\in H^1 (X, T_X)$. We say that the deformations of $X$ are unobstructed if all the obstructions to deformations are $0$. 

If the deformations of $X$ are unobstructed (i.e., $\dim T_0 \Def (X) = \dim \Def (X)$), then the base $\Def (X)$ is a small open neighborhood of the origin in $H^1 (X,T_X)$. The following theorem is due to Bogomolov in the hyperk\"ahler case and to Tian-Todorov in the general case.

\begin{theorem}
If the canonical bundle $K_X$ is trivial (we say $X$ is $K$-trivial), then the deformations of $X$ are unobstructed.
\end{theorem}

We have the following facts.

\begin{itemize}
\item If $X$ is K\"ahler, then so is any small deformation of $X$.

\item If $X$ is K\"ahler and $K$-trivial, then small deformations $X_t$ of $X$ are also K\"ahler and $K$-trivial and $h^1 (T_{X_t})$ is constant.

\item If $X$ is holomorphic symplectic, then small deformations of $X$ are also holomorphic symplectic. If $X$ is irreducible holomorphic symplectic, then all fibers of any deformation of $X$ are irreducible holomorphic symplectic.
\end{itemize}

\subsection{The Beauville-Bogomolov form}

The key to understanding the deformations of hyperk\"ahler manifolds is the period domain. Small open subsets of the period domain are isomorphic to $\Def (X)$. We define the period domain using the second cohomology of hyperk\"ahler manifolds, together with a non-degenerate quadratic form: the Beauville-Bogomolov form.

Suppose $X$ is irreducible holomorphic symplectic (irreducible hyperk\"ahler) of dimension $2n$ and choose $\sigma \in H^0 (\Omega_X^2)$ such that
\[
\int_X (\sigma \osigma )^n = 1.
\]
For $\alpha \in H^2 (X, \bC)$, define
\[
q_X (\alpha ) := \frac{n}{2} \int_X \alpha^2 (\sigma \osigma )^{n-1} + (1-n) \int_X \sigma^{n-1}\osigma^n \alpha \int_X \sigma^n \osigma^{n-1} \oalpha.
\]
One can show this is equal to
\[
q_X (\alpha ) = \lambda \mu + \frac{n}{2} \int_X \beta^2 (\sigma \osigma )^{n-1}
\]
where $\alpha = \lambda \sigma + \beta + \mu \osigma$ with $\beta \in H^{1,1} (X)$.

Beauville showed that there exists $d_X \in \bN$ such that
\[
\int_X \alpha^{2n} = d_X (q_X (\alpha))^n.
\]
Therefore, if $r_X$ is the positive real root of $d_X$, then $\tq_X := r_X q_X$ is an $n$-th root of the $n$-th power cup-product on $H^2(X, \bC)$.

The quadratic form $\tq_X$ is integer valued on $H^2 (X, \bZ)$, indivisible, non-degenerate, of signature $(3, b_2 - 3)$ on $H^2 (X, \bR)$. Furthermore,
\[
\tq_X (\sigma ) =0, \quad \tq_X (\sigma + \osigma ) >0
\]
and
\[
\tq_X (\sigma_t ) =0, \quad \tq_X (\sigma_t + \osigma_t ) >0
\]
for $t$ close to $0$ in any deformation of $X$.

The form $\tq$ is called the Beauville-Bogomolov form of the hyperk\"ahler manifold. The inspiration for the Beauville-Bogomolov form came from the study of the Fano variety of lines of a cubic fourfold. There, it naturally appears as the intersection form on the fourth cohomology of the cubic threefold which is isomorphic to the second cohomology of its Fano variety of lines which is a hyperk\"ahler manifold.

Note that for $n=1$, $\tq_X = 2q_X$ is the usual intersection form on $H^2 (X, \bZ)$.

\subsection{The local period domain and local Torelli}
Define
\[
Q_X := \{ \alpha \mid q_X (\alpha ) =0 , q_X (\alpha + \oalpha ) >0 \} \subset \oQ_X \subset \bP H^2 (X, \bC).
\]
We saw that for $t\in \Def (X)$ close to $0$, $q_X (\sigma_t ) =0 , q_X (\sigma_t + \osigma_t ) >0$. Hence we can define the local period map
\[
\begin{array}{rcl}
P_X : \Def (X) & \lra & Q_X \\
t & \longmapsto &[\sigma_t].
\end{array}
\]
This is holomorphic because $\langle \sigma_t\rangle = H^{2,0} (X_t) = H^0 (\Omega^2_{X_t})$ varies holomorphically with $t$: $H^0 (\Omega^2_{X_t})$ is the fiber of the holomorphic line bundle $\pi_*\Omega^2_{\cX / \Def (X)}$ on $\Def (X)$.

We have the local Torelli theorem:
\begin{theorem}
The local Torelli map $P_X$ is a local isomorphism, i.e., $d P_X$ is an isomorphism at $0$.
\end{theorem}

\subsection{The period domain}

We now construct the global period domain for hyperk\"ahler manifolds. For this we first fix the discrete data of a lattice which will usually be abstractly isomorphic to the second integral cohomology of a hyperk\"ahler manifold with its Beauville-Bogomolov form.

\begin{definition}
A lattice is the data of a free $\bZ$-module $\Gamma$ of finite rank with an integral non-degenerate quadratic form $q_\Gamma$.
\end{definition}

\begin{definition}
Given a lattice $(\Gamma, q_\Gamma)$, the period domain $Q_\Gamma$ is
\[
Q_\Gamma := \{ \alpha \mid q_X (\alpha ) =0 , q_X (\alpha + \oalpha ) >0 \} \subset \oQ_\Gamma \subset \bP (\Gamma \otimes_\bZ \bC).
\]
\end{definition}

\subsection{The moduli space of marked holomorphic symplectic manifolds and local period maps}

We will construct a moduli space of marked holomorphic symplectic manifolds and a global period map on it which is, roughly speaking, a glueing of local period maps.

\begin{definition}
\begin{enumerate}
\item A marking of an irreducible holomorphic symplectic manifold is a lattice isomorphism
\[
\varphi : (H^2 (X, \bZ), \tq_X ) \stackrel{\cong}{\lra} (\Gamma, q_\Gamma).
\]
\item The pair $(X, \varphi)$ is called a marked manifold.
\item Two marked manifolds $(X, \varphi)$, $(X', \varphi')$ are isomorphic if there exists $f : X \ra X'$ such that $\varphi' = \varphi \circ f^*$. We write $(X, \varphi) \cong (X', \varphi')$.
\item The moduli space of marked irreducible holomorphic symplectic manifolds is the set
\[
\cM_\Gamma := \{ (X, \varphi ) \} /\cong.
\]
\end{enumerate}
\end{definition}
We use the local period map to show that $\cM_\Gamma$ is a smooth (non-Hausdorff) complex analytic space:

Given an irreducible holomorphic manifold $X$, choose a marking $\varphi : H^2 (X, \bZ) \ra \Gamma$.
The Kuranishi family $\cX \ra \Def (X)$ is locally isomorphic to the period domain $Q_\Gamma$: The marking $\varphi : H^2 (X, \bZ) \ra \Gamma$ induces isomorphisms forming the commutative diagram
\[
\xymatrix{
         Q_X \ar@{^{(}->}[d]_{} \ar[r]^{\cong} & Q_\Gamma \ar@{^{(}->}[d]\\
         \bP H^2 (X, \bC) \ar[r]^{\cong}      & \bP (\Gamma \otimes_\bZ \bC).}
\]
Hence an open ball in the Kuranishi space $\Def (X)$ is isomorphic to an open ball in $Q_\Gamma$. Such open balls cover $\cM_\Gamma$ and the analytic structures on intersections coincide because the Kuranishi family is the local universal deformation of all of its fibers. Hence we obtain a well-defined smooth complex analytic structure on $\cM_\Gamma$.

\subsection{The global period map and Verbitsky's global Torelli theorem}

\begin{definition}
The global period map is
\[
\begin{array}{cccc}
         P : & \cM_\Gamma & \lra & Q_\Gamma \subset \oQ_\Gamma \subset \bP (\Gamma \otimes_\bZ \bC)\\
         & (X, \varphi) &  \longmapsto    & [\varphi (\sigma)].
         \end{array}
\]
\end{definition}
Verbitsky's global Torelli theorem \cite{Verbitsky2013} (also see \cite{Huybrechts2012tor} and \cite{Looijenga2021teich}) for compact hyperk\"ahler manifolds is the following.
\begin{theorem}
The map $P$ is generically injective on each connected component of $\cM_\Gamma$.
\end{theorem}
Note that the datum of the line $H^{2,0} (X) \subset H^2 (X, \bC)$ determines the Hodge structure on $H^2(X, \bZ)$: $H^{0,2}(X) = \overline{H^{2,0}(X)}$ (complex conjugate), $H^{2,0}(X)^\perp = H^{2,0} (X) \oplus H^{1,1}(X)$, $H^{1,1}(X) = \left( H^{2,0} (X) \oplus H^{1,1}(X) \right) \cap \overline{\left( H^{2,0} (X) \oplus H^{1,1}(X) \right)}$.

We say that the global Torelli theorem holds for a class of manifolds, if a manifold is determined by its Hodge structure, possibly together with the data of a polarization (such as the form $\tq_X$ in the hyperk\"ahler case). For instance, two complex tori are isomorphic if and only if their first cohomologies are isomorphic as Hodge structures. Two Riemann surfaces are isomorphic if and only if their first cohomologies are Hodge isometric, i.e., they are isomorphic as Hodge structures and, under the given Hodge isomorphism, the intersection forms for the two curves coincide. Similarly, two K3 surfaces are isomorphic if their second cohomologies are Hodge isometric.

In fact we have stronger Torelli theorems in the above cases: for complex tori, any Hodge isomorphism between the first cohomologies of two tori is {\em induced} by an isomorphism of the tori. For curves, any Hodge isometry between their first cohomologies is induced by an isomorphism between the curves {\em up to a change of sign}. For {\em generic} K3 surfaces, any Hodge isometry between the second cohomologies is induced by an isomorphism of the surfaces up to a sign.

For hyperk\"ahler manifolds of dimension $>4$, none of the above stronger versions of Torelli hold. There are examples of
\begin{enumerate}
\item non-isomorphic (but bimeromorphic) compact hyperk\"ahler manifolds with Hodge isometric second cohomologies \cite{Debarre1984},
\item non-birational projective hyperk\"ahler manifolds of dimension $4$ with Hodge isometric second cohomologies, \cite{Namikawa2002}.
\end{enumerate}
\begin{question}
Is there a good characterization of irreducible holomorphic symplectic manifolds that are Hodge isometric but not isomorphic?
\end{question}
We have the following maps of moduli spaces
\[
\xymatrix{
        \Teich (X) \ar@{->}[d]_{} \ar@{=}[r] & \{ \hbox{complex structures on } X \} /\sim^0 & \\
         \cM_\Gamma (X) \ar[d] \ar@{=}[r] &  \{ \hbox{marked complex structures on } X \} /\approx & \\
        \cM_{cx} (X) \ar@{=}[r] &  \{ \hbox{complex structures on } X \} /\sim \ar@{=}[r] & \Teich (X) / G }
\]
and the period map
\[
\xymatrixcolsep{5pc}\xymatrix{
\Teich (X) \ar[r]^{local\; isom.} & \cM_\Gamma(X) \ar[r]^{P_\Gamma} & \quad Q_\Gamma & \mkern-72mu\subset \oQ_\Gamma \subset \bP (\Gamma \otimes \bC).}
\]
The spaces $\Teich (X)$ and $\cM_\Gamma (X)$ are non Hausdorff smooth analytic spaces and $Q_\Gamma$ is a (Hausdorff) simply connected complex manifold. Verbitsky constructed a new (Hausdorff) complex manifold $\cM_\Gamma^s (X)$ which is obtained by identifying all non-separated points of $\cM_\Gamma (X)$. In other words
\[
\cM_\Gamma^s (X) = \cM_\Gamma (X) / \equiv
\]
where, for two points $p,q \in \cM_\Gamma (X)$, $p\equiv q$ when every neighborhood of $p$ contains $q$ and every neighborhood of $q$ contains $p$.
The period map then factors through $\cM_\Gamma^s (X)$:
\[
\xymatrixcolsep{5pc}\xymatrix{
P_\Gamma : \quad \cM_\Gamma (X) \ar[r]^{local\; isom.} & \cM_\Gamma^s(X) \ar[r]^{P_\Gamma^s} & \quad Q_\Gamma.}
\]
Verbitsky proved
\begin{theorem}
The map $P^s_\Gamma$ is surjective from any connected component of $\cM_\Gamma^s(X)$ to $Q_\Gamma$.
\end{theorem}
Combined with the facts that $P^s_\Gamma$ is a local isomorphism and $Q_\Gamma$ is simply connected, this implies
\begin{corollary}
The map $P^s_\Gamma$ induces an isomorphism from any connected component of $\cM_\Gamma^s(X)$ to $Q_\Gamma$.
\end{corollary}
Verbitsky's proof uses twistor conics which we will describe in the next section.

The following results of Huybrechts help us understand the difference between $\cM_\Gamma(X)$ and $\cM_\Gamma^s(X)$.
\begin{proposition}
If two marked hyperk\"ahler manifolds $(X, \varphi)$ and $(X', \varphi')$ correspond to two non-separated points of $\cM_\Gamma(X)$, then $X$ and $Y$ are bimeromorphic and their period $P_\Gamma (X,\varphi) = P_\Gamma (X',\varphi)$ is contained in the hyperplane $Q_\Gamma \cap \alpha^\perp$ for some $\alpha \in \Gamma$.
\end{proposition}
\begin{proposition}
Suppose given a bimeromorphism $f : X \ra X'$ between compact, hyperk\"ahler manifolds. Then there exists families of compact hyperk\"ahler manifolds
\[
\cX \lra D, \quad \quad \cX' \lra D
\]
over a complex disc $D$ such that
\begin{enumerate}
\item $\cX_0 \cong X$ and $\cX'_0 \cong X'$,
\item there exists a bimeromorphism $F : \cX \ra \cX'$ commuting with the projections to $D$ which is an isomorphism over $D\setminus \{0\}$ and induces $f$ on $\cX_0 \cong X \dashrightarrow \cX'_0 \cong X'$.
\end{enumerate}
\end{proposition}
\begin{proposition}
For any $x\in Q_\Gamma$, the set of hyperk\"ahler complex structures on a differentiable manifold $X$ with period $x\in Q_\Gamma$ consists of a finite number of bimeromorphic equivalence classes.
\end{proposition}

\section{Twistor spaces and twistor conics}

\subsection{Hyperk\"ahler structures}\label{subsectwistor}
Given $X$ hyperk\"ahler, let $g$ be the hyperk\"ahler metric of $X$. We saw that there exists complex structures $I, J, K$ such that $g$ is K\"ahler with respect to $I,J,K$ and $IJK = -1$. In fact $g$ is K\"ahler with respect to any linear combination $\lambda = aI + bJ + cK$ such that $a^2 + b^2 + c^2 =1$. The K\"ahler form associated to $\lambda$ is $\omega_\lambda ( \cdot, \cdot ) := g (\lambda \cdot, \cdot)$. So we have a family $\{ (X, \lambda) \mid \lambda \in S^2 \}$ of compact K\"ahler manifolds.

\subsection{Twistor spaces}

With the notation above, the twistor space $\cX \ra \bP^1$ of $(X, g)$ is the product $X\times \bP^1$ (as a real manifold) endowed with the almost complex structure
\[
\begin{array}{cccc}
         I_{X \times \bP^1} : & T_x X \oplus T_\lambda \bP^1 & \lra & T_x X \oplus T_\lambda \bP^1 \\
         & (v, w) &  \longmapsto    & (\lambda (v), I_{\bP^1} (w))
         \end{array}
\]
which is integrable by a result of Hitchin, Karlhede, Lindstr\"om, Ro{\v c}ek.

\subsection{Twistor conics}

Fix a lattice $(\Gamma, q_\Gamma)$, isometric to $(H^2 (X, \bZ), \tq_X)$. Recall that the signature of $q_\Gamma \otimes \bR$ is $(3, b_2 -3)$ where $b_2$ is the second Betti number of $X$. Since $\bP^1$ is simply connected, we can choose consistent markings on all the fibers of $\cX \ra \bP^1$ to obtain the period map
\[
\begin{array}{cccc}
P_g : & \bP^1 & \lra & Q_\Gamma \\
 & \lambda & \longmapsto & [\sigma_{(X, \lambda )}]
\end{array}
\]
whose image is a twistor conic.

One can show that it is the intersection of a linearly embedded $P = \bP^2$ with $Q_\Gamma$ in $\bP (\Gamma \otimes \bC)$. Furthermore $P = \bP (W\otimes \bC)$ where $W$ is a three dimensional real subspace of $\Gamma \otimes \bR$ totally positive for the intersection form $q_\Gamma$.

Conversely, one can show that each choice of a $3$-dimensional real space $W \subset \Gamma\otimes \bR$ positive for $q_\Gamma$ gives a twistor conic:
\[
C := \bP (W\otimes \bC) \cap Q_\Gamma \subset Q_\gamma.
\]
Recall the following
\begin{definition}
A K\"ahler class is the cohomology class of a $(1,1)$ form which is K\"ahler with respect to some metric. The K\"ahler cone is the cone generated by all K\"ahler classes.
\end{definition}
A consequence of the Calabi-Yau theorem is the following.
\begin{corollary}
Suppose $(M,I,g)$ is compact H\"ahler with $c_1 (K_M) =0$. Then, in each K\"ahler class on $M$, there exists a unique Ricci-flat metric. Furthermore, the Ricci-flat K\"ahler metrics on $M$ form a smooth family of dimension $h^{1,1} (M)$ isomorphic to the K\"ahler cone.
\end{corollary}
Therefore, given the family $\{ (X, \lambda) \mid \lambda \in S^2 \}$ as in \ref{subsectwistor}, 
for every K\"ahler class $\alpha \in H^{1,1}(M)$, there exists a unique hyperl\:ahler metric $g_\lambda$, K\"ahler with respect to $\lambda$, such that $[\omega_{g_\lambda}] = \alpha$.

For each such metric $[\omega_{g_\lambda}]$, we can construct a twistor family. In other words, through each point $[(X,I)]$ of the twistor conic there passes another twistor conic.

One can show

\begin{proposition}
$Q_\Gamma$ is twistor path connected, i.e., any two points of $Q_\Gamma$ can be joined by a connected sequence of twistor conics.
\end{proposition}
From which it follows
\begin{corollary}
The period map $P_\Gamma :\cM_\Gamma \ra Q_\Gamma$ is surjective on any connected component of $\cM_\Gamma$.
\end{corollary}

\subsection{Hyperholomorphic bundles}

We start with the definition of hyperholomorphic bundles.
\begin{definition}
Given a hermitian vector bundle $B$ on $X$, with hermitian connection $\theta$, we say $(B,\theta)$ is hyperholomorphic if it is compatible with all the complex structures $\lambda\in S^2 = \bP^1$.
\end{definition}
\begin{definition}
A $C^\infty$ vector bundle $B$ on $X$ is hermitian if it has a hermitian metric (denoted $\langle, \rangle$). A connection
\[
\theta : B \lra B \otimes T_X^*
\]
is hermitian if the metric is (covariantly) constant with respect to $\theta$. If we are given a complex structure $I$ on $B$, we say that $\theta$ and $I$ are compatible if the curvature form
\[
\Theta : B \lra B \otimes \Lambda^2 T_M^*
\]
is a $(1,1)$-form with respect to $I$.
\end{definition}
Intuitively, considering the twistor family
\[
\xymatrix{
X \times \bP^1 \ar@{=}[r]_{C^\infty} & \cX \ar[d] &  \ar[l]^{C^\infty} B\times \bP^1 \\
 & \bP^1,
}
\]
the $C^\infty$ vector bundle $B\times \bP^1$ on $\cX$ has a structure of complex vector bundle holomorphic on each fiber $(X, \lambda)$ of $\cX \ra \bP^1$.

Stability conditions allow us to construct moduli spaces of bundles.
\begin{definition}
Fix a K\"ahler form $\omega$ on $X$. For a coherent sheaf $F$ on $X$, put
\[
\deg (F) := \frac{1}{\vol(X)} \int_X c_1 (F) \wedge \omega^{n-1}
\]
where $n$ is the complex dimension of $X$ and $\vol(X) := \int_X \omega^n$. Define
\[
\slope (F) := \frac{\deg (F)}{\rank(F)}
\]
where $\rank (F)$ is the complex rank of $F$.
We say $F$ is stable with respect to $\omega$ if for all subsheaves $F' \subset F$ with $\rank (F') < \rank (F)$, we have
\[
\slope (F' ) < \slope (F).
\]
We say $F$ is semi-stable with respect to $\omega$ if for all subsheaves $F' \subset F$, we have
\[
\slope (F' ) \leq \slope (F).
\]
\end{definition}

Verbitsky (see \cite{VerbitskyKaledin1999}) proved that, given a vector bundle $B$ on $(X,I)$, if $c_1 (B)$ and $c_2(B)$ are of type $(1,1)$ and $(2,2)$ with respect to all complex structures $\lambda\in S^2 =\bP^1$ on $X$, then $B$ is hyperholomorphic. In particular, the class $c_2(B)$ is analytic on each $(X, \lambda)$.

A useful characterization of stable bundles is given by the Hitchin-Kobayashi correspondence. To state it, we first need the following definition.
\begin{definition}
Let $\omega$ be the K\"ahler form of $M$ and denote by $\Lambda : \Omega^{1,1}_M \otimes B \ra B$ the adjoint of cup-product with $\omega$. A hermitian metric with curvature form $\Theta : B \ra B \otimes \Omega^{1,1}_M$ is Hermitian-Einstein if the composition $\Lambda \Theta : B \ra B$ is a multiple of the identity.
\end{definition}
The Hitchin-Kobayashi correspondence, proved by Donaldson, Uhlenbeck and Yau is the following theorem.
\begin{theorem}
Suppose $B$ is an indecomposable bundle on a compact K\"ahler manifold $M$. Then $B$ is stable if and only if $B$ has a Hermitian-Einstein metric.
\end{theorem}

\section{Examples of hyperk\"ahlers in dimension $2$ and beyond, by Samir Canning}

\subsection{Betti and Hodge numbers of $K3$ surfaces}
The purpose of this exercise is to compute the Betti and Hodge numbers of a complex $K3$ surface $X$, which is the simplest example of a hyperk\"ahler manifold. 
Feel free to add the additional assumption that $X$ is algebraic if you are more comfortable in that setting. 

\begin{problem}
Show that $H^0(X,\bZ)=H^4(X,\bZ)=\bZ$, $H^1(X,\bZ)=0$, and $H^3(X,\bZ)$ is torsion. (Hint: use the exponential exact sequence.)
\end{problem}
\begin{problem}
Show that $H^2(X,\bZ)$ is torsion free. Conclude that $H^3(X,\bZ)=0$. (Hint: continue analyzing the exponential exact sequence, using that $\Pic(X)$ is torsion free. Prove this if you know about Riemann-Roch. For the second statement, use the universal coefficient theorem for cohomology.)
\end{problem}

Recall the Hirzebruch--Riemann--Roch Theorem.
\begin{theorem}[Hirzebruch--Riemann--Roch]
Let $F$ be a (holomorphic) vector bundle on a compact complex manifold $X$. Then, 
\[
\chi(X,F)=\int_X \ch(F)\td(X).
\]
\end{theorem}
When we write $c_i(X)$, we mean $c_i(T_X)$, where $T_X$ is the tangent bundle. Here are the first few terms of the Chern character and Todd class for reference:
\[
\ch(F)=\rank(F)+c_1(F)+\frac{1}{2}(c_1(F)^2-2c_2(F))+\cdots
\]
and
\[
\td(F)=1+\frac{1}{2}c_1(F)^2+\frac{1}{12}(c_1(F)^2+c_2(F))+\cdots
\]
\begin{problem}
Compute $c_2(X)$ for $X$ a $K3$ surface. (Hint: set $F=\cO_X$.)
\end{problem}

\begin{problem}
Compute $H^2(X,\bZ)$. (Hint: take $F=\Omega_X$.)
\end{problem}

You have now computed all of the Betti numbers. Next, we will compute the Hodge numbers.  
\begin{definition}
Let $X$ be a compact K\"ahler manifold. The Hodge numbers of $X$ are 
\[
h^{p,q}=\dim H^q(X,\Omega^p_X).
\]
\end{definition}

\begin{theorem}[The Hodge Decomposition]
Let $X$ be a compact K\"ahler manifold. There is a direct sum decomposition
\[
H^i(X,\bZ)\otimes \mathbb{C}=H^i(X,\mathbb{C})=\bigoplus_{p+q=i} H^q(X,\Omega^p_X).
\]
Moreover $h^{p,q}=h^{q,p}$.
\end{theorem}
\begin{problem}
Compute all of the Hodge numbers of a compact complex $K3$ surface $X$. 
\end{problem}

\begin{further}
The same ideas, especially the use of the Hirzebruch--Riemann--Roch Theorem, can be used to give restrictions on the Betti and Hodge numbers of higher dimensional hyperk\"ahler manifolds. For more in this direction, see the paper of Salamon \cite{Salamon} and Debarre's exposition thereof \cite{Debarre}. For even further restrictions on the Betti numbers of hyperk\"ahler fourfolds, see the paper of Guan \cite{Guan}.
\end{further}

\subsection{Identifying hyperk\"ahler manifolds}
One of the most interesting areas of research in hyperk\"ahler geometry is the construction of examples. This exercise will focus on identifying examples. We begin with some basic problems.
\begin{problem}
Convince yourself that any holomorphic two-form $\sigma$ on a complex manifold $X$ induces a morphism of bundles
\[
\sigma: T_X\rightarrow \Omega^1_X.
\]
where $T_X$ is the tangent bundle and $\Omega^1_X$ is the cotangent bundle.
\end{problem}
We call $\sigma$ non-degenerate if the morphism above is an isomorphism.

\begin{problem}
Can you convince yourself that K3 surfaces are irreducible hyperk\"ahler? (Hint: the tricky part is probably the simply connectedness. It may require some extra background knowledge.)
\end{problem}

\begin{problem}
Show that $h^{2,0}=h^{0,2}=1$, $K_X\cong \cO_X$, and that $\dim(X)$ is even for any irreducible compact hyperk\"ahler manifold $X$. 
\end{problem}

Now that we know that $K_X$ is trivial for compact hyperk\"ahler manifolds $X$, a natural question is: given a $K_X$-trivial manifold, how can we show that it is hyperk\"ahler, if it is? We will focus on a real-life example due to Debarre--Voisin \cite{DebarreVoisin}. The same type of argument works for another famous example of Beauville--Donagi \cite{BeauvilleDonagi} (the Fano variety of lines on a cubic fourfold.)

Let $V_{10}$ be a $10$-dimensional complex vector space. Let $\omega\in \wedge^3 V_{10}^\vee$ be a $3$-form on $V_{10}$. We define a subvariety of $G(6,V_{10})$:
\[
X_{\omega}:=\{[W]\in G(6,V_{10}): \omega|_{W\times W\times W}\equiv 0\}.
\]
\begin{problem}
Show that for a general choice of $\omega$, $X_{\omega}$ is a smooth fourfold. (Hint: show that $X_{\omega}$ is given by the vanishing of a section of a certain globally generated vector bundle.)
\end{problem}

\begin{problem}
Show that $K_{X_{\omega}}\cong \cO_{X_{\omega}}$. (Hint: use adjunction.) 
\end{problem}
Now that we know we have a $K_X$-trivial variety, we want to show it's hyperk\"ahler. 
Using something called the Koszul resolution, one can compute the Euler characteristic of the structure sheaf:
\[
\chi(X_{\omega},\cO_{X_{\omega}})=3.
\]
\begin{definition}
A strict Calabi-Yau manifold is a simply connected projective manifold $X$ such that $H^0(X,\Omega^p_X)=0$ for $0<p<\dim(X)$.
\end{definition}

\begin{problem}
Show that any simply connected smooth $K_X$-trivial compact K\"ahler fourfold with $\chi(X,\cO_X)=3$ is irreducible compact hyperk\"ahler. (Hint: use the nice multiplicative properties of $\chi(X,\cO_X)$.)
\end{problem}

\begin{further}
The proof that $X_{\omega}$ above is hyperk\"ahler is done differently (more geometrically) in \cite{DebarreVoisin}. I also highly recommend the classic paper \cite{BeauvilleDonagi}. It turns out in both cases, the resulting hyperk\"ahler is deformation equivalent to the Hilbert scheme of $2$ points on a $K3$ surface. 
\end{further}

\section{Basic properties of Lagrangian fibrations of Hyperk\"ahlers, by Yajnaseni Dutta}

The following exercises are based on a couple of fundamental results from \cite{Matsushita1999fib} and \cite{Matsushita2005lagr}. Given a Lagrangian fibration $f\colon X\to B$ of a Hyperk\"ahler manifold $X$, the geometry and topology of $B$ are heavily influenced by $X$. In fact, Matsushita conjectured that $B\simeq \bP^n$. It is known by work of Hwang \cite{Hwa08} that if $B$ is smooth then $B\simeq \bP^n$. The conjecture is known to be true if $\dim B = 2$ by recent results of \cite{BK18, HX20, Ou19}

\subsection{Lagrangian fibrations}
Let $S$ be a K3 surface and $f\colon S\to C$ a proper surjective morphism on to a smooth irreducible curve with connected fibres \footnote{we will call such a morphism \textsl{fibration} throughout the rest of these exercises.}. 
		\begin{problem} Show that $C\simeq \bP^1$. \hint{Use that $S$ is simply connected.}
		\end{problem}
	
	\begin{problem} Show that the general fibres of $f$ are elliptic curves. \hint{Use Adjunction.}
		\end{problem}
	
\begin{problem} Find an explicit fibration of the Fermat quartic $(x^4+y^4+z^4+w^4=0) \subset \bP^3$. \hint{rewrite as equality of two fractions.}\end{problem}

\vspace*{.5em} Let $X$ be a hyperk\"ahler manifold of dimension $2n$. The following exercises show how similar the situation is in higher dimensions. The quadratic space $(H^2(X,\bR), q_X)$ controls much of the geometry of $X$ and is a central gadget in the study of hyperk\"ahler manifolds. 

Recall that $q_X$ is a priori dependant on the symplectic form $\sigma\in H^0(X,\Omega_X^2)$, however, up to scaling, it is independent of $\sigma$. Here are some key properties of $q_X$ (we denote the associated bilinear form again by $q_X$).
\begin{itemize}
	\item The normalized symplectic form $\sigma$ satisfies $q_X(\sigma) = 0$ and $q_X(\sigma+\overline{\sigma}) =1$.	
	\item More generally, for 
	$\alpha_i\in H^2(X)$, we have \[\displaystyle\int_X\alpha_1\cdots\alpha_{2n} = c_X\sum_{s\in S_n}q_X(\alpha_{s(1)}, \alpha_{s(2)})\dots q_X(\alpha_{s(2n-1)}, \alpha_{s(2n-2)})\] for some constant $c_X$ depending only on $X$. As a consequence, we obtain $\int_X\sigma\overline{\sigma}\omega^{2n-2} = c'q_X(\omega)^{n-1}$.
	\item If a line bundle $L$ is ample, then $q_X(c_1(L)) > 0$. The K\"ahler cone is contained in a connected component of $\{\alpha\in H^{1,1}(X,\bR)\mid q_X(\alpha)>0\}$. Partial converses to these statements exist. For instance, if $L$ is a line bundle with $q_X(L)>0$ then $X$ is projective \cite[Prop.\ 26.13]{GrossHuybrechtsJoyce2003}. Furthermore, if $q_X(\alpha) > 0$ and, for every rational curve $C\subset X$, $\int_C \alpha > 0$, then $\alpha$ is a K\"ahler class \cite[Th\'eor\`eme 1.2]{Bou01}. 
	\item $H^{1,1}(X,\bC)$ is orthogonal to $H^{2,0}(X,\bC) \oplus H^{0,2}(X,\bC)$ with respect to $q_X$.
	\item By \cite{Bog96, Ver96} whenever there exists $0\neq \beta\in H^2(X,\bC)$ that satisfies $q_X(\beta) = 0$, we have $\beta^n\neq 0$ and $\beta^{n+1} = 0$
\end{itemize}
\noindent We begin with a Hodge index type theoerem.

\begin{problem}\label{ex:zerointersection} Given a divisor $E$ on $X$, show that if $E$ satisfies $E^{2n} = 0$ and $E\cdot A^{2n-1} = 0$ for some ample bundle $A$, then $E\sim 0$. (Hint: Use $q_X(tE+A) = t^2q(E) + 2tq(E,A)+q(A)$ for any variable $t$ and that $(tE+A)^{2n} = c_X q_X(tE+A)^n$.)
\end{problem}
	
	\begin{problem}\label{ex:intersection} Given a divisor $E$ on $X$, show that if $E$ satisfies $E^{2n} = 0$ and $E\cdot A^{2n-1}>0$ for some ample line bundle $A$, then $q_X(E,A)>0$ and the following are true
	\[\begin{split}
	E^m\cdot A^{2n-m} &= 0 \text{ ; for } m>n\\
	E^m\cdot A^{2n-m} &> 0 \text{ ; for } m\leq n.
	\end{split}
	\]
(Hint: Expand $q_X(tE+A) $ as in the previous exercise.)
	\end{problem}
	
	\begin{problem} Let $f\colon X\to B$ be a fibration of a hyperk\"ahler manifold $X$\footnote{you may assume both $X$ and $B$ are projective, although the results presented here work in a more general setting.}.
Using the previous exercise show that $\dim B = n$. (Hint: Apply the previous exercise to the pull-back of an ample class $H$ on $B$.)
		\end{problem}

	\begin{problem}\label{ex:ns} Show that $\Pic(B)$ is of rank 1. (Hint: Show that any divisor $E$ on $X$ that satisfies $E^{2n} = 0$ and $E^n\cdot (f^*H)^n = 0$ is in fact a rational multiple of $f^*H$.)
	\end{problem}

		\noindent For the next exercise we need the definition of a Lagrangian (possibly singular) subvariety. Recall that
	\begin{definition}
	A subvariety $Y\subset X$ is said to be a Lagrangian subvariety if $\dim Y=\frac{1}{2}\dim X$ and there exists a resolution of singularities $\mu\colon Y'\to Y$ such that $\mu^*\sigma|_Y = 0$.
	\end{definition}

	\begin{problem} Show that a general fibre of $f$ is Lagrangian. By a classical theorem, the general fibres of $f$ are then complex tori. A more recent result of Voisin \cite[Prop.\ 2.1]{Cam06} or, more generally, \cite[Theorem 1.1]{Leh16}, shows that even if $X$ is not projective, a Lagrangian subvariety of a hyperk\"ahler manifold is always projective. Thus, a general fibre $F$ is isomorphic to an abelian variety. \hint{Let $A$ be an ample class on $X$. Argue that if $\sigma|_F\neq 0$, then $\int_F(\sigma\overline{\sigma}|_F)\cdot A|_F^{n-2}\neq 0$. Using the polynomial $q_X(t f^*H+A)^{n-1}$ and the second property listed above show that it is in fact 0.}
	\end{problem}

	\begin{problem} Show that \textsl{every} fibre of $f$ is Lagrangian and hence $f$ is equidimensional. (Hint: Use the map $H^2(X, \mathcal{O}_X)\to H^0(B, R^2f_* \mathcal{O} _X)$ induced by the Leray spectral sequence and that $R^2f_* \mathcal{O} _X$ is torsion free.)
	\end{problem}

	\begin{problem} Show that $B$ is $\bQ$-factorial with at worst Kawamata log terminal singularities. \hint{Use that $f$ is equidimensional and \cite[Lemma 5.16]{KM98} which states that if the source of a finite surjective map between normal varieties is $\bQ$-factorial and klt then so is the target.}
	\end{problem}

	\vspace*{0.5em}
	\noindent For the next exercise, recall and use the following
	\begin{definition}[Kodaira Dimension]\label{def:kd}
Let $X$ be a $\bQ$-factorial variety. Then 
\[\kappa(X)=\sup_{m} \dim\overline{\phi_m(X)}\]
where $\phi_m:X\dashrightarrow \bP^{P_m}$ is the rational map defined 
by the global sections of $\omega_X^{\otimes m}$
and $P_m=\dim H^0(X, \omega_X^{\otimes m})$.
Another way to interpret this is \[\kappa(X)\coloneqq\text{trdeg}_{_k}\left (\bigoplus_mH^0(X,\omega_X^{\otimes m})\right)-1\]
 where the algebra structure on the right side is given by the multiplication map. 
\end{definition}
Iitaka's $C_{n,m}$ conjecture then states that 
\begin{conjecture}
Let $f:X\to B$ be a fibration of smooth projective varieties of dimension $n$ and $m$, respectively, and let $F$ be a general fibre of $f$. Then,
\[\kappa(X)\ge \kappa(F) + \kappa(B).\]
\label{thm:iitaka}
\end{conjecture}
\noindent By a result of Kawamata \cite[Theorem 1.1(2)]{Kaw85}, the conjecture is known when $F$ is a minimal variety.

\begin{problem} Assume $B$ is smooth, show that $B$ is Fano, i.e., the inverse of the canonical bundle of $B$ is ample. (Hint: use that the Picard rank of $B$ is 1 and Kawamata's result above.)
	\end{problem}
	
	\begin{problem} Assume $B$ is smooth. Let $B^0$ be the open set where $f$ is smooth. Let $X^0 \coloneqq f^{-1}(B^0)$. Show that $R^if^0_*\cO_{X^0} = \Omega_{B^0}^i$. (Hint: Use $\Omega_{X^0}^1\simeq \mathcal{T} _{X^0}$ to conclude that $f^* \mathcal{T} _{B^0}\simeq \Omega_{X^0/B^0}^1$.)
		\end{problem}
		Matsushita \cite{Matsushita2005lagr} (also see \cite{Matsushita1999fib}) extends this equality to the big open set $U$ which includes the smooth points of the discriminant divisor $D_f$, using Deligne's canonical extension. Then, using the reflexivity of $R^if_*\cO_X$ and the isomorphism $R^nf_*\cO_X\simeq \omega_B$, he shows that $R^if_*\cO_X\simeq \Omega_B^i$.

\section{Rational curves on K3 surfaces and Euler characteristics of Moduli spaces, by David Stapleton}

We work through an idea of Beauville \cite{Beauville1999rat}, following work of Yau and Zaslow \cite{YauZaslow}, which uses hyperk\"ahler geometry to count the number of rational curves in a very general K3 surface of degree 2d.\\

\noindent{\textbf{Problem 1.}} Assume that a K3 surface $X$ admits an \textit{elliptic pencil} -- that is a map
\[
\pi\cl X\ra \bP^1
\]
so that the general fibers are smooth genus 1 curves. Assume that all the fibers that do not have geometric genus 1 are irreducible rational curves with a single node. Count the number of rational fibers. (Hint: If $R = \sqcup_{i=1}^n R_i$ is the union of rational curves, compute the topological Euler characteristic using the formula:
\[
e(X) = e(R) + e(X\setminus R)
\]
and compute $e(R_i)$.)

\subsection{Hyperk\"ahlers as moduli spaces of sheaves on K3 surfaces.} Let $X$ be a very general K3 surface of degree 2d with primitive line bundle $L$ (with $L^2 = 2d$) and let $\Pi = \bP(H^0(X,L))\cong \bP^{d-1}$. Moduli spaces of sheaves on $X$ are frequently hyperk\"ahler manifolds. Here are two examples:
\begin{enumerate}
\item\textbf{Hilbert schemes of $n$ points} on $X$ -- denoted $X^{[n]}$, this space compactifies the space of unordered distinct points on $X$ by considering length $n$ subschemes as their limits.
\item\textbf{Compactified Jacobians} -- denoted $\ocJ^d(X)$ -- parametrizing coherent sheaves supported on curves $C\in \Pi$, which when thought of as sheaves on $C$ are line bundles (or torsion-free sheaves of rank 1 when $C$ is singular) of degree $d$.
\end{enumerate}

\vspace{12pt}
\noindent{\textbf{Problem 3.}} Show that if $X$ is a K3 surface, then $\Pi$ contains only finitely many rational curves (curves with geometric genus 0).

\vspace{12pt}
\noindent{\textbf{Problem 4.}} Compute the dimension of $X^{[n]}$ and $\ocJ^d(X)$.

\vspace{12pt}
\noindent{\textbf{Problem 5.}} Show that the hyperk\"ahlers $X^{[g]}$ and $\ocJ^g(X)$ are birational.

\vspace{12pt}

There is a natural map
\[
\pi\cl \ocJ^g(X)\ra \Pi
\]
which sends a coherent sheaf $\cF$ to the curve in $\Pi$ that it is supported on.

\vspace{12pt}

\noindent{\textbf{Problem 6.}} Show that the general fiber of $\pi$ is an Abelian variety. Describe the fibers over a general point $C\in \Pi$.

\vspace{12pt}

\noindent{\textbf{Problem 7.}} (this is \cite[Prop. 2.2]{Beauville1999rat}) Let $C$ be an integral curve such that the normalization $\hC$ has genus $\ge 1$. We show that $e(\ocJ^d(C))=0$ as follows.\\
(1) Find a line bundle $\cM$ on $C$ which is torsion of order $m$ (for any $m>0$). (This uses the comparison between the Jacobian of $C$ and of $\hC$.)\\
(2) Show that tensoring by $\cM$ is a free action of $\bZ/m\bZ$ on $\ocJ^d(\hC)$.\\
(3) Conclude that $m$ divides $e(\ocJ^d(C))$ for all $m>0$.

\vspace{12pt}

It follows by the scissor property of Euler characteristics that
\[
e(\ocJ^g(X)) =\sum_{R_i \in \Pi}e(\ocJ^g(R_i))
\]
where $R_i\in \Pi$ is a rational curve and $\pi^{-1}(R_i)$ is the fiber over $R_i$ (i.e., the set of torsion free sheaves of rank 1 and degree $g$ supported on $R_i$).\\

\noindent{\textbf{Problem 8.}} Show that
\[
e(\ocJ^g(R_i))=1
\]
if $R_i$ is a nodal, irreducible rational curve. (Thus by a result of Xi Chen \cite{XiChen}, if $X$ is very general then
\[
e(\ocJ^g(X)) = \#\{ R_i\in\Pi\}.)
\]
Hint: Locally at a node $p\in R_i$ there are only 2 types of rank 1 torsion free sheaves (1) line bundles and (2) the ideal sheaf of a point. Show that if $p_1,\cdots,p_g\in R_i$ are the nodes then $\ocJ^g(R_i)$ is stratified into loci $\ocJ^g_S\subset \ocJ^g(R_i)$ consisting of torsion-free sheaves that are not locally free exactly at the points in a subset $S\subset \{p_1,\cdots,p_g\}$. Conclude that the only stratum where $e(\ocJ^g_S)\ne 0$ is when $S=\{p_1,\cdots, p_g\}$ (a single point). See also \cite[\S 3]{Beauville1999rat}.

\vspace{12pt}

It remains to actually calculate the Euler characteristic of $\ocJ^g(X)$. This relies on
\begin{enumerate}
\item The birational invariance of Euler characteristic for hyperk\"ahlers (see \cite{HuybrechtsDeformation} or use the birational invariance of betti numbers of Calabi-Yaus \cite{batyrev}).
\item The computation of the Euler characteristic of $X^{[n]}$ by G\"ottsche \cite{Gottsche} (see \cite{DeCataldo2000hilb} for a nice explanation of these results).
\end{enumerate}

In particular, for a K3 surface, by (1) and (2) we have:
\begin{center}
$\sum (\#\text{ rational curves on a K3 of genus $g$})q^g= \sum_{g\ge 0} e(\ocJ^g(X))q^g$\\
$=\sum_{g\ge 0} e(X^{[g]})q^g = \Pi_{k=1}^\infty \left(\frac{1}{1-q^k}\right)^{e(X)}$
\end{center}
\noindent where the sum over $g\ge 0$ is understood to take a very general K3 surface of genus $g$.\\

\noindent{\textbf{Problem 9.}} Compute the Euler characteristic of $X^{[2]}$ for any complex surface using that

(1) there is a birational map
\[
h: X^{[2]}\ra X^{(2)}
\]
to the symmetric product $X^{(2)}:= X^2/\Sigma_2$ which is given by blowing up the diagonal locus and

(2) the exceptional divisor of $h$ is a $\bP^1$-bundle over $X$.\\

\noindent{\textbf{Problem 10.}} Find the number of bitangents to a very general plane sextic curve $C\subset \bP^2$ using that a very general K3 surface of genus $2$ is a double cover of $\bP^2$ branched at such a sextic.


\providecommand{\bysame}{\leavevmode\hbox to3em{\hrulefill}\thinspace}
\providecommand{\MR}{\relax\ifhmode\unskip\space\fi MR }
\providecommand{\MRhref}[2]{%
  \href{http://www.ams.org/mathscinet-getitem?mr=#1}{#2}
}
\providecommand{\href}[2]{#2}

\end{document}